\documentclass[12pt, oneside,reqno]{amsart}	
\usepackage{geometry}                		
\usepackage{graphicx}				
\usepackage{amssymb}
\usepackage{amsmath}
\usepackage[all]{xy}
\usepackage{mathtools}
\usepackage{tikz-cd}
\usepackage{enumitem}
\usepackage{tikz-cd}

\newtheorem{thm}{Theorem}[section]
\newtheorem{lem}[thm]{Lemma}
\newtheorem{cor}[thm]{Corollary}
\newtheorem{prop}[thm]{Proposition}

\newtheorem{defn}{Definition}[section]

\numberwithin{equation}{section}

\newcommand{\Balpha}{\mbox{$\hspace{0.12em}\shortmid\hspace{-0.62em}\alpha$}} 

\def\Pb{\ifmmode{\Bbb P}\else{$\Bbb P$}\fi}
\def\Z{\ifmmode{\Bbb Z}\else{$\Bbb Z$}\fi}
\def\C{\ifmmode{\Bbb C}\else{$\Bbb C$}\fi}
\def\R{\ifmmode{\Bbb R}\else{$\Bbb R$}\fi}
\def\S{\ifmmode{S^2}\else{$S^2$}\fi}

\def\S{\cal S}

\newenvironment{pf}{\paragraph{Proof:}}{\hfill$\square$ \newline}

\begin {document}
	
\title{Low Entropy and the Mean Curvature Flow with Surgery}

\begin{abstract}In this article, we extend the mean curvature flow with surgery to mean convex hypersurfaces with entropy less than $\Lambda_{n-2}$. In particular, 2-convexity is not assumed. Next we show the surgery flow with just the initial convexity assumption $H - \frac{\langle x, \nu \rangle}{2} > 0$ is possible and as an application we use the surgery flow to show that smooth $n$-dimensional closed self shrinkers with entropy less than $\Lambda_{n-2}$ are isotopic to the round $n$-sphere.  \end{abstract}

\author {Alexander Mramor \& Shengwen Wang}
\address{Department of Mathematics, University of California Irvine, Irvine, CA 92617}
\address{Department of Mathematics, SUNY Binghamton, Binghamton, NY 13902}
\email{mramora@uci.edu,swang@math.binghamton.edu}

\maketitle

\section{Introduction.}
In \cite{CM} Colding and Minicozzi introduced a quantity $\lambda(M)$ of a submanifold $M \subset \R^{N}$ they called the entropy; its an especially interesting quantity for a number of reasons, one of which being that it in a strong sense captures information about a submanifold, the intuition being somehow low entropy surfaces should be simpler. On the other hand the set of submanifolds is a robust set under perturbations compared to an apriori curvature condition, which constrain the geometry of a submanifold pointwise. We denote by $\Lambda_k=\lambda(\mathbb S^k)=\lambda(\mathbb S^k\times\mathbb R^{n-k})$. According to Stone's computation \cite{St}: $$\Lambda_1>\frac{3}{2}>\Lambda_2>...>\Lambda_n\rightarrow\sqrt2$$.

For a given dimension $n$, it is already known that $\Lambda_n$ is a lower bound on entropy for a hypersurface $M^n \subset \R^{n+1}$, that this is obtained exactly when $M$ is a round sphere, and that surfaces with entropy close to $\Lambda_n$ are Hausdorff close to a round sphere - see \cite{CIMW,BW,BW0,SW}.  In $\mathbb R^4$, from work due to Bernstein and Lu Wang \cite{BW1} when the entropy is below $\Lambda_2=\lambda(\mathbb S^2\times\mathbb R)$ any closed hypersurface $M$ with $\lambda(M) < \Lambda_2$ (so entropy between $\Lambda_2$ and $\Lambda_3$) is topologically a 3-sphere, and later Bernstein and the second named author showed in \cite{BSW} the level-set flow such surfaces with entropy less than $\Lambda_2$ stay connected until extinction. It is a natural question to ask what can be said of flows of surfaces with the next highest level of entropy, namely those with entropy below $\Lambda_1$ or, more generally, $\Lambda_{n-2}$.
$\medskip$

As a step towards answering this question, note that mean convex self shrinkers of entropy bounded by $\Lambda_{n-2}$ are either $S^n$ or $S^{n-1} \times \R$ and hence are 2-convex. In addition, mean curvature flow with surgery has been established for globally $2$-convex hypersurfaces. This suggests a mean curvature flow with surgery for mean convex low entropy (but not necessarily 2-convex) hypersurfaces is possible and is the topic of this article. In this article we confirm this and show how to extend the mean curvature flow with surgery as defined by Haslhofer and Kleiner to mean convex surfaces of low entropy:
\begin{thm} Let $\mathcal{M} = \mathcal{M}(\alpha, \gamma, n, \Pi)$ be the set of $\alpha$-noncollapsed closed hypersurfaces in $\R^{n+1}$ with entropy less than $\Pi < \Lambda_{n-2}$ and $H < \gamma$, then for any $M \in \mathcal{M}$ there is a mean curvature flow with surgery for a uniform choice of parameters $H_{th}, H_{neck}, H_{trig}$. 
\end{thm}
Note that if $\Pi < \Lambda_{n-1}$ then any $M \in \mathcal{M}$ must shrink to a point; we will assume throughout that $\Pi \gneq \lambda_{n-1}$. In the process of proving our theorem we construct a mean convex hypersurface of low entropy (which in this article refers to surfaces $M$ with $\lambda(M) < \Lambda_{n-2}$ unless otherwise stated) that develops a neckpinch, showing that surgeries are to be expected for $M \in \mathcal{M}$ above; see section 3.2 below. There are several corollaries of the surgery; the first consequence was first noted in the original mean curvature flow with surgery paper by Huisken and Sinestrari \cite{HS2}:
\begin{cor} If $M \in \mathcal{M}$, then $M \cong S^n$ or a finite connect sum of $S^{n-1} \times S^1$ \end{cor}
The second corollary, an extension of the first corollary, was observed in the 2-convex case by the first named author in \cite{Mra}; only mean convexity of the MCF was required so the proof immediately adapts to the low entropy case:
\begin{cor} Let $\Sigma(d, C, \mathcal{M})$ be the set of hypersurfaces $ M \in \mathcal{M}$ with diam$(M) < d$ (or equivalently up to translation, $M  \subset B_d(0)$) and $H < C$. Then $\Sigma(d, C, \mathcal{M})$ up to isotopy consists of finitely many hypersurfaces, in fact at most $2^{2n} \frac{( 12dC\sqrt{n})^n}{\alpha^n}$. 
\end{cor}
We turn next to applications of the flow to self shrinkers of low entropy. Perturbing a low entropy self shrinker as in \cite{CIMW} we obtain a low entropy surface with $H - \frac{\langle x, \nu \rangle}{2} > 0$. From \cite{Smk,Lin} a related convexity assumption is preserved under the flow and in fact adapting some estimates of Lin \cite{Lin} we may then run the surgery essentially using that high curvature regions will still be mean convex. Combining it with a recent result of Hershkovits and White \cite{HW} yields:
\begin{thm} Suppose $M^n$ is a closed self shrinker with entropy bounded above by $\Lambda_{n-2}$, where $n \geq 3$. Then $M$ is diffeomorphic to $S^n$ and is in fact isotopic to round $S^n$. 
\end{thm}
So we see that in some sense the self shrinkers of entropy bounded by $\Lambda_{n-2}$ have properties one expects from general hypersurfaces of entropy bounded by $\Lambda_{n-1}$. 
$\medskip$

Foundational results in this direction concerning self shrinkers include the paper \cite{CIMW} of Colding, Minicozzi, Ilmanen, and White, where the self shrinking sphere is first classified as the closed self shrinker of lowest entropy (in fact, our result is an improvement on theorem 0.6 in their paper). A more recent result is due to Hershkovits and White \cite{HW}, where a rigidity theorem for self shrinkers relating entropy and topological type is given and which plays an important role in the proof of our application - their result in fact says that a self shrinker with the above entropy bound is weakly homotopic to a sphere. 
$\medskip$

The structure of the article is as follows: first we give preliminary information on the mean curvature flow with surgery and Colding and Minicozzi's entropy, then we describe the proof of theorem 1.1 above. Next the application to low entropy self shrinkers is described. We end with some concluding remarks about future avenues of investigation and a supplementary appendix justifying the roundabout construction of the ``low entropy'' surgery cap.  
$\medskip$

$\textbf{Acknowledgements:}$ The authors are indebted to Jacob Bernstein, the second named author's advisor, who encouraged them to investigate applications of their flow with surgery to self shrinkers of low entropy. They also thank Or Hershkovits for pointing out a mistake in an earlier draft and his helpful comments concerning theorem 1.4, and Longzhi Lin for discussions concerning his paper. The first named author additionally thanks his advisor, Richard Schoen, for his support. Finally they thank the referee for their careful reading and insightful questions which helped improve the exposition of this article. 
 
\section{Preliminaries.} 
The first subsection introducing the mean curvature flow we borrow quite liberally from the first named author's previous paper \cite{Mra} although a small introduction to self shrinkers has also been included. The second subsection concerns the mean curvature flow with surgery as constructed by Haslhofer and Kleiner in \cite{HK1}. The third subsection introduces some basic facts and definitions concerning Colding and Minicozzi's entropy introduced in \cite{CM}. 
\subsection{Classical formulation of the mean curvature flow}
In this subsection we start with the differential geometric, or ``classical," definition of mean curvature flow for smooth embedded hypersurfaces of $\R^{n+1}$; for a nice introduction, see \cite{Mant}. Let $M$ be an $n$ dimensional manifold and let $F: M \to \R^{n+1}$ be an embedding of $M$ realizing it as a smooth closed hypersurface of Euclidean space - which by abuse of notation we also refer to $M$. Then the mean curvature flow of $M$ is given by $\hat{F}: M \times [0,T) \to \R^{n+1}$ satisfying (where $\nu$ is outward pointing normal and $H$ is the mean curvature):
\begin{equation}
\frac{d\hat{F}}{dt} = -H \nu, \text{ } \hat{F}(M, 0) = F(M)
\end{equation} 
(It follows from the Jordan separation theorem that closed embedded hypersurfaces are oriented). Denote $\hat{F}(\cdot, t) = \hat{F}_t$, and further denote by $M_t$ the image of $\hat{F}_t$ (so $M_0 = M$). It turns out that (2.1) is a degenerate parabolic system of equations so take some work to show short term existence (to see its degenerate, any tangential perturbation of $F$ is a mean curvature flow). More specifically, where $g$ is the induced metric on $M$:
\begin{equation}
 \Delta_g F = g^{ij}(\frac{\partial^2 F}{\partial x^i \partial x^j} - \Gamma_{ij}^k \frac{\partial F}{\partial x^k}) = g^{ij} h_{ij} \nu = H\nu
\end{equation}
Now one could apply for example deTurck's trick to reduce the problem to a nondegenerate parabolic PDE (see for example chapter 3 of \cite{Bake}) or similarly reduce the problem to an easier PDE by writing $M$ as a graph over a reference manifold by Huisken and Polden (see \cite{Mant}). At any rate, we have short term existence for compact manifolds.
$\medskip$

Now that we have established existence of the flow in cases important to us, let's record associated evolution equations for some of the usual geometric quantities: 
\begin{enumerate}
\item $\frac{\partial}{\partial t} g_{ij} = - 2H h_{ij}$
$\medskip$

\item $\frac{\partial}{\partial t} d\mu = -H^2 d\mu$
$\medskip$

\item $\frac{\partial}{\partial t} h^i_j = \Delta h^i_j + |A|^2 h^i_j$
$\medskip$

\item $\frac{\partial}{\partial t} H = \Delta H + |A|^2 H$
$\medskip$

\item $\frac{\partial}{\partial t} |A|^2 = \Delta |A|^2 - 2|\nabla A|^2  + 2|A|^4$
\end{enumerate}
So, for example, from the heat equation for $H$ one sees by the maximum principle that if $H > 0$ initially it remains so under the flow. There is also a more complicated tensor maximum principle by Hamilton originally developed for the Ricci flow (see \cite{Ham1}) that says essentially that if $M$ is a compact manifold one has the following evolution equation for a tensor $S$: 
\begin{equation}
\frac{\partial S}{\partial t} = \Delta S + \Phi(S)
\end{equation} 
and if $S$ belongs to a convex cone of tensors, then if solutions to the system of ODE
\begin{equation}
\frac{\partial S}{\partial t} = \Phi(S)
\end{equation}
stay in that cone then solutions to the PDE (2.2) stay in the cone too (essentially this is because $\Delta$ ``averages"). So, for example, one can see then that convex surfaces stay convex under the flow very easily this way using the evolution equation above for the Weingarten operator. Similarly one can see that $\textbf{2-convex hypersurface}$ (i.e. for the two smallest principal curvatures $\lambda_1, \lambda_2$, $\lambda_1 + \lambda_2 > \beta H$ everywhere for some $\beta > 0$) remain 2-convex under the flow.
$\medskip$

Another important curvature condition in this paper is $\alpha \textbf{ non-collapsing}$: a mean convex hypersurface $M$ is said to be 2-sided $\alpha$ non-collapsed for some $\alpha > 0$ if at every point $p\in M$, there is an interior and exterior ball of radius $\alpha/H(p)$ touching $M$ precisely at $p$. This condition is used in the formulation of the finiteness theorem. It was shown by Ben Andrews in \cite{BA} to be preserved under the flow for compact surfaces. (a sharp version of this statement, first shown by Brendle in \cite{Binsc} and later Haslhofer and Kleiner in \cite{HK3}, is important in \cite{BH} where MCF+surgery to $n=2$ was first accomplished). 
$\medskip$

Finally, perhaps the most geometric manifestation of the maximum principle is that if two compact hypersurfaces are disjoint initially they remain so under the flow; this fact is used in section 4.2 below. So, by putting a large hypersphere around $M$ and noting under the mean curvature flow that such a sphere collapses to a point in finite time, the flow of $M$ must not be defined past a certain time either in that as $t \to T$, $M_t$ converge to a set that isn't a manifold.  
$\medskip$

Note this implies as $t \to T$ that $|A|^2 \to \infty$ at a sequence of points on $M_t$; if not then we could use curvature bounds to attain a smooth limit $M_T$ which we can then flow further, contradicting our choice of $T$. Thus weak solutions to the flow are necessitated; one type of weak solution is the mean curvature flow with surgery explained below. 
$\medskip$

This also implies that understanding the singularities to the mean curvature flow is a topic of great importance in the field. Self shrinkers, solutions to the elliptic equation
\begin{equation}
H = \frac{\langle x, \nu \rangle}{2}
\end{equation} 
Are in a very strong sense the singularity models for the flow; namely when one does a tangent flow blowup at a point, the $t = -1$ time slice of the limit flow is modeled on a solution to (2.3) by Huisken's monotonicity formula (see \cite{H}). Considering the Gaussian metric $g_{ij} = e^{\frac{-|x|^2}{4}} \delta_{ij}$ on $\R^{n+1}$, one easily sees by calculating the first variation formula for a hypersurface in this metric that self shrinkers are precisely the minimal surfaces. This metric shares some similarities with Ricci positive spaces, for example that minimal surfaces (self shrinkers) in this metric satisfy a Frenkel property (that is, any two compact minimal surface must intersect, see \cite{F}) as exploited by the authors in the previous article \cite{MW} (the Frankel property for noncompact self shrinkers has also been explored, see \cite{IPM} -- although their theorem is not used one can also prove some unknottedness theorems for noncompact self shrinkers as discussed in \cite{Mra2}). 
$\medskip$

Its known by work of Huisken and later Colding and Minicozzi than under very mild assumptions the only mean convex self shrinkers are generalized round cylinders $S^k \times \R^{n-k}$, but outside of the mean convex setting little is understood. There are a patchwork of non mean convex examples, such as those constructed by ODE methods \cite{Ang}, minmax methods \cite{Ket}, and gluing methods \cite{KKM}. Theorem 1.4 says however that all compact examples with topology must have ``high" entropy.

\subsection{Mean curvature flow with surgery for compact 2-convex hypersurfaces in $\R^{n+1}$}
First we give the definition of $\Balpha$ controlled:
\begin{defn}(Definition 1.15 in \cite{HK1}) Let $\Balpha = (\alpha, \beta, \gamma) \in (0, N-2) \times (0, \frac{1}{N-2}) \times (0, \infty)$. A smooth compact closed subamnifold $M^n \subset \R^{n+1}$ is said to be $\Balpha$-controlled if it satisfies
\begin{enumerate}
\item is $\alpha$-noncollapsed
\item  $\lambda_1 + \lambda_2 \geq \beta H$ ($\beta$ 2-convex)
\item $H \leq \gamma$
\end{enumerate}
\end{defn} 
Speaking very roughly, for the mean curvature flow with surgery approach of Haslhofer and Kleiner, like with the Huisken and Sinestrari approach there are three main constants, $H_{th} \leq H_{neck} \leq H_{trig}$. If $H_{trig}$ is reached somewhere during the mean curvature flow $M_t$ of a manifold $M$ it turns out the nearby regions will be ``neck-like" and one can cut and glue in appropriate caps (maintaining 2-convexity, etc) so that after the surgery the result has mean curvature bounded by $H_{th}$. The high curvature regions have well understood geometry and are discarded and the mean curvature flow with surgery proceeds starting from the low curvature leftovers. Before stating a more precise statement we are forced to introduce a couple more definitions. First an abbreviated definition of the most general type of piecewise smooth flow we will consider. 
\begin{defn} (see Definition 1.3 in \cite{HK1}) An $(\alpha, \delta)-\textit{flow}$ $M_t$ is a collection of finitely smooth $\alpha$-noncollapsed flows $\{M_t^i \cap  U\}_{t \in [t_{i-1}, t_i]}$, $(i = 1, \ldots k$; $t_0 < \ldots t_k$) in an open set $U \subset \R^{n+1}$, such that:
\begin{enumerate}
\item for each $i = 1, \ldots, k-1$, the final time slices of some collection of disjoint strong $\delta$-necks (see below) are replaced by standard caps, giving $M^{\#}_{t_i} \subset M^i_{t_i} =: M_{t_i}^-$ (in terms of the regions they bound). 
$\medskip$

\item the initial time slice of the next flow, $M_{t_i}^{i+1} =: M_{t_i}^{+}$, is obtained from $M_{t_i}^\#$ by discarding some connected components.
\end{enumerate} 
\end{defn}
Of course, now we should define what we mean by standard caps, cutting and pasting, and strong $\delta$-necks. Since we will need them in the sequel, we will give the full definitions; these are essentially definitions 2.2 through 2.4 in \cite{HK1}:
\begin{defn} A standard cap is a smooth convex domain that coincides with a smooth round half-cylinder of radius 1 outside a ball of radius 10. 
\end{defn}
The model we give for a standard cap will morally agree with the definition given above although the radius outside which it will agree with the round cylinder will potentially need to be taken larger than 10. In the next definition note that in practice (considering a neck point $p$ on $M$) $s$ will be $\frac{n-1}{H(p)}$; in particular after rescaling it will be equal to n-1: 
\begin{defn}  We say than an $(\alpha, \delta)$-flow $M_t$ has a strong $\delta$-neck with center $p$ and radius $s$ at time $t_0 \in I$, if $\{s^{-1} \cdot (M_{t_0 + s^2 t} - p) \}_{t \in (-1, 0]}$ is $\delta$-close in $C^{[1/\delta]}$ in $B_{1/\delta}^U \times (-1, 0]$ to the evolution of a round cylinder $S^{n} \times \R$ with radius 1 at $t = 0$, where $B^U_{1/\delta} = s^{-1} \cdot ((B(p, s/\delta) \cap U) - p) \subset B(0, 1/\delta) \subset \R^{n+1}$. 
\end{defn}  
Now is the definition of cutting and pasting:
\begin{defn} We say that a final time slice of a strong $\delta$-neck ($\delta \leq \frac{1}{10 \Gamma}$) with center $p$ and radius $s$ is replaced by a pair of standard caps if the pre-surgery domain $M^-$ is replaced by a post surgery domain $M^+$  such that 
\begin{enumerate}
\item the modification takes place inside a ball $B = B(p, 5 \Gamma s)$
\item there are bounds for the second fundamental form and its derivatives:
\begin{center}
$\sup\limits_{M^+ \cap B} |\nabla ^\ell A| \leq C_\ell s^{-1 - \ell}$   ($\ell = 0, 1, 2, \ldots$)
\end{center}
\item if $B \subset U$, then for every point $p_+ \in M^+ \cap B$ with $\lambda_1(p_+) < 0$, there is a point $p_\# \in M^\# \cap B$ with $\frac{\lambda_1}{H}(p_+) < \frac{\lambda_1}{H}(p_\#)$
\item if $B(p, 10\Gamma s) \subset U$, then $s^{-1}(M^+ - p)$ is $\delta'(\delta)$-close in n$B(0, 10\Gamma)$ to a pair of disjoint standard saps that are at distance $\Gamma$ from the origin. 
\end{enumerate}
\end{defn}

With these definitions in mind before moving on we state an important set of properities that standard caps satisfy. As long as the cap we construct satisfies the defintion of standard cap above and that after the gluing the postgluing domain adheres to definition 2.5 above the proposition will be true. We include this though for completeness sake since it is used, as one may check, many times in the proof of the canonical neighborhood theorem.
\begin{prop} Let $C$ be a standard cap with $\alpha, \beta > 0$. There is a unique mean curvature flow $\{C_t\}_{t \in [0, 1/2(N-2))}$ starting at $C$. It has the following properties.
\begin{enumerate}
\item It is $\alpha$-noncollapsed, convex, and $\beta$-uniformly 2-convex.
\item There are continuous increasing functions $\underline H, \overline H: [0, \frac{1}{2(N-2)} \to \R$, with $H(t) \to \infty$ as $t \to \frac{1}{2(N-2)}$ such that $\underline H(t) \leq H(p,t) \leq \overline{H}(t)$ for all $p \in C_t$ and $t \in [0, 1/2(N-2))$.
\item For every $\epsilon > 0$ and $\tau < \frac{1}{2(N-2)}$ there exists an $R = R(\epsilon, \tau) < \infty$ such that outside $B(0,R)$ the flow $C_t$, $t \in [0, \tau]$, is $\epsilon$ close the the flow of the round cylinder. 
\item For every $\epsilon > 0$, there exists a $\tau = \tau(\epsilon) < \frac{1}{2(N-2)}$ such that every point $(p,t) \in \partial K_t$ with $t \geq \tau$ is $\epsilon$-close to a $\beta$-uniformly 2-convex ancient $\alpha$-noncollapsed flow. 
\end{enumerate}
\end{prop}
We sketch the proof of canonical neighborhood theorem below (of course, full details are in \cite{HK1}). Before that we finally state the main existence result of Haslhofer and Kleiner; see theorem 1.21 in \cite{HK1}
\begin{thm} (Existence of mean curvature flow with surgery). There are constants $\overline{\delta} = \overline{\delta}(\Balpha) > 0$ and $\Theta(\delta) = \Theta(\Balpha, \delta) < \infty$ ($\delta \leq \overline{\delta}$) with the following significance. If $\delta \leq \overline{\delta}$ and $\mathbb{H} = (H_{trig}, H_{neck}, H_{th})$ are positive numbers with $H_{trig}/H_{neck}, H_{neck}/H_{th}, H_{neck} \geq \Theta(\delta)$, then there exists an $(\Balpha, \delta, \mathbb{H})$-flow $\{M_t\}_{t \in [0, \infty)}$ for every $\Balpha$-controlled surface $M$.  
\end{thm} 
The reason we choose to employ the scheme set out by Haslhofer and Kleiner because the surgery problem is then reduced to showing ancient $\alpha$-nonocollapsed flows of suitably low entropy are $\beta$ 2-convex for some $\beta > 0$. Without going into more details than necessary, we recall one last theorem we will need in the sequel, see theorem 1.22 in \cite{HK1}:
\begin{thm} (Canonical neighborhood theorem) For all $\epsilon > 0$, there exists $\overline{\delta} = \overline{\delta}(\Balpha) > 0$, $H_{can}(\epsilon) = H_{can}(\Balpha, \epsilon) < \infty$ and $\Theta_\epsilon(\delta) = \Theta_\epsilon(\Balpha, \delta) < \infty$ ($\delta \leq \overline{\delta}$) with the following signifigance. If $\delta < \overline{\delta}$ and $M$ is an $(\Balpha, \delta, \mathbb{H})$-flow with $H_{trig}/H_{neck}, H_{neck}/H_{th} \geq \Theta_\epsilon(\delta)$, then any $(p,t) \in \delta M$ with $H(p,t) \geq H_{can}(\epsilon)$ is $\epsilon$-close to either (a) a $\beta$-uniformly 2-convex ancient $\alpha$-noncollapsed flow, or (b) the evolution of a standard cap preceeded by the evolution of a round cylinder. 
\end{thm}
$\medskip$

The above theorem is roughly proven by letting the surgery ratios above degenerate to infinity for a sequence of flows and analyizing the possibilities for the limits, which are guaranteed by a convergence theorem of Haslhofer and Kleiner. In the case of no surgeries the limit that is ancient and $\beta$ two convex and $\alpha$ non collapsed, so that the theorem follows since the convergence is in a suitably strong topology. If there are surgeries, then it follows that the limit contains a line (more specifically, see claim 4.3 and the discussion afterwards in \cite{HK1}), from which $(b)$ follows. this part uses the properties of the cap that are satisfied in proposition 2.1 above. 
$\medskip$

The last case does not employ two convexity, so to see that the canonical neighborhood theorem is true in our setting it suffices to show that ancient, $\alpha$-noncollapsed, low entropy flows are in fact $\beta$ 2-convex for some $\beta > 0$ and that our cap is suitably constructed to satisfy proposition 2.1; these are both attended to in the next section. 
$\medskip$

Now, to prove the existence of the surgery, Haslhofer and Kleiner proceed by finding regions which seperate high curvature regions, where some points have $H = H_{trig}$, and low curvature regions where $H \leq H_{th}$; see claim 4.6 in \cite{HK1}. These will contain strong neck points in the sense above on which they can do surgery; see claim 4.7 in \cite{HK1}. 
$\medskip$

 If the ancient flow found in the canonical neighborhood theorem is compact, it will be diffeomorphic to a sphere, see the discussion after claim 4.8 in \cite{HK1}. Furthermore as long as $H_{th}$ is taken large enough (roughly large enough to employ the canonical neighborhood theorem for appropriately small $\epsilon$ as we do in section 3.3 below) all points in the intermediate region between $H_{th}$ and $H_{neck}$ can be forced to be neck points; we will also refer to this region as the neck region below. This is essentially also contained in the argument in the proof of corollary 1.25 in \cite{HK1} following claim 4.8 therein. 
 $\medskip$
 
 For readers perhaps more familiar with the approach to surgery of Huisken and Sinestrari in \cite{HS2}, this is essentially the content of their neck continuation theorem (more precisely, theorem 8.1 in \cite{HS2}); one starts by finding a neck point, and the statement is essentially that one may continue the neck as long as $H$ is large (in our context, $H > H_{th}$), $\lambda_1/H$ is small, and there are no previous surgeries in the way. If the second or third conditions are violated the case then is that the neck is ended by a convex cap. 
\subsection{Background on Colding and Minicozzi's entropy.}

In \cite{CM} Colding and Minicozzi discovered a useful new quantity called the entropy to study the mean curvature flow. To elaborate, consider a hypersurface $\Sigma^k \subset \R^{\ell}$; then given $x_0 \in \R^{\ell}$ and $r > 0$ define the functional $F_{x_0, r}$ by 
\begin{equation}
F_{x_0, r}(\Sigma) = \frac{1}{(4 \pi r)^{k/2}} \int_\Sigma e^\frac{-|x - x_0|^2}{4r} d\mu
\end{equation} 
Colding and Minicozzi then define the entropy $\lambda(\Sigma)$ of a submanifold to be the supremum over all $F_{x_0, r}$ functionals:
\begin{equation} 
\lambda(\Sigma) = \sup\limits_{x_0, r} F_{x_0, r}(\Sigma) 
\end{equation} 
Important for below is to note that equivalently $\lambda(\Sigma)$ is the supremum of $F_{0,1}$ when we vary over rescalings (changing $r$) and translations (choice of $x_0$). For hypersurfaces with polynomial growth this supremum is attained and, for self shrinkers $\Sigma$, $\lambda(\Sigma) = F_{0,1}(\Sigma)$. In fact, self shrinkers are critical points for the entropy so it is natural next to ask what the stable ones are. If $\Sigma_2$ is a normal variation of $\Sigma$ and $x_s, t_s$ are variations with $x_0=0, r = 1$, 
\begin{equation}
\partial_s \mid_{s = 0} \Sigma_s = f\nu, \partial_s \mid_{s=0} x_s = y, \text{ and }\partial_s \mid_{s = 0} t_s = h
\end{equation} 
The second variation formula one find is:
\begin{equation}
 F_{0,1} = (4 \pi)^{-n/2} \int_\Sigma (-fLf + 2fhH - h^2H^2 f \langle y, \nu \rangle - \frac{\langle y, \nu \rangle^2}{2})e^{\frac{-|x|^2}{4}} d\mu
 \end{equation}
 where $L$ is given by the following: 
\begin{equation}
L = \Delta + |A|^2  -\frac{1}{2} \langle x, \nabla(\cdot) \rangle + \frac{1}{2}
\end{equation} 
One can easily check that, where $v$ is a constant vector field on $\R^n$, both $\langle v, \nu \rangle$ and $H$ are eigenfunctions with eigenvalues $-1, -\frac{1}{2}$ respectively for $L$; $LH = H$ and $L\langle v, \nu \rangle = \frac{1}{2}\langle v, \nu \rangle$. $L$ is self adjoint in the weighted space $L^2(e^{\frac{-|x|^2}{4}})$, so has a discrete set of eigenvalues with corresponding orthogonal sets of eigenfunctions. If a self shrinker isn't mean convex $H$ switches signs on $\Sigma$, so by the minmax characterization for eigenvalues on a surface $\Sigma$ must not be the lowest eigenvalue, and that there is a positive function $f$ that is $L^2(e^{\frac{-|x|^2}{4}})$ orthogonal to both $H$ and $\langle v, \nu \rangle$. It follows from the second variation formula that $f$ gives rise to an entropy decreasing variation of $\Sigma$, so that namely $\Sigma$ is not stable. Thus all stable self shrinkers (with some area growth assumptions) are mean convex and must be spheres and cylinders; more precisely: 
\begin{thm}(Theorem 0.12 in \cite{CM}) Suppose that $\Sigma$ is a smooth complete embedded self-shrinker without boundary and with polynomial volume growth.
\begin{enumerate} 
\item If $\Sigma$ is not equal to $S^k \times \R^{n-k}$, then there is a graph $\widetilde{\Sigma}$ over $\Sigma$ of a function with arbitrarily small $C^m$ norm (for any fixed $m$) so that $\lambda(\widetilde{\Sigma}) < \lambda(\Sigma)$
\item If $\Sigma$ is not $S^n$ and does not split off a line, then the function in (1) can be taken to have compact support. 
\end{enumerate}
\end{thm} 
 This theorem has been extended to the singular settng by Zhu in \cite{Z} (see theorem 0.2). Furthermore the entropy is monotone decreasing under the flow by Huisken monotonicity \cite{H} so, if the entropy of a surface is lower than that of a certain self shrinker, that self shrinker won't be the singularity model for any singularities of the surface under the flow later on - this of course is essential and a surgery flow for low entropy surfaces wouldn't be sensible otherwise. 
 $\medskip$
 
 We end this discussion with a lemma which restricts which $F$-functionals we will need to consider when estimating the entropy. This is contained in the argument of lemma 7.7 of \cite{CM} which says that the entropy is achieved by an $F$ functional for a smooth closed embedded hypersurface. 
 \begin{lem}\label{conhull} Let $\Sigma \subset \R^{n+1}$ be smooth and embedded. For a given $r > 0$, the supremum over $x_0$ of $F_{x_0, r}(\Sigma)$ is achieved within the convex hull of $\Sigma$ 
 \end{lem}
 
To see this, Colding and Minicozzi note that from the first variation $F_{x_0, r}$ must be a critical point for fixed $r$ when the integral $x - x_0$ vanishes, which couldn't occur if $x_0$ wasn't in the convex hull of $\Sigma$.

 \section{Proof of Theorem 1.1.}
The proof of theorem 1.1 amounts to showing the following two things, most of the work in the article being to establish (2): 
\begin{enumerate}
\item Ancient mean convex solutions of low entropy are in fact uniformly $\beta$ 2-convex for some $\beta > 0$, and
\item The low entropy condition is preserved across surgeries. 
\end{enumerate}
Using the first item one can proceed exactly as in \cite{HK1} to establish the canonical neighborhood theorem and so on as discussed in section 2.2. More precisely since all the low entropy, mean convex, ancient hypersurfaces are uniformly 2-convex, all the statements in section 3 of \cite{HK1} are true in our setting. 

In the following (this concerns the second step) without loss of generality we will assume there is only one surgery performed at a time for a time slice $T$; if there are multiple to be performed at once the argument below works if they are considered successively (within a fixed time slice).

\subsection{Structure of $\alpha$-noncollapsed ancient flows of low entropy}
First we quickly establish item (1) above. Before proceeding we recall that $\alpha$-noncollapsing, entropy, and $\beta$ 2-convexity are all scale invariant conditions/quantities. We apply the next proposition with $\epsilon_0 = \Lambda_{n-2} - \Pi$:
\begin{prop} Pick $\epsilon_0 > 0$. There exists $\beta>0$, depending only on $\alpha,n$, and $\epsilon_0$ such that if $M^n_t$ be an $\alpha$-noncollapsed ancient flow in $\mathbb R^{n+1}$ and $\lambda(M_{t})<\Lambda_{n-2} - \epsilon_0$ then there exists some $\beta$ so that $M_t$ is $\beta$ 2-convex.
\end{prop}
\begin{proof}
Suppose not, there exists a sequence of ancient $\alpha$-noncollapsed flows $\{M_{i,t}\}$ with $(p_i,t_i)\in\{M_{i,t}\}$ such that $\frac{\lambda_1(p_i)+\lambda_2(p_i)}{H}<\frac{1}{i}\rightarrow0$. We can translate and rescale to get a sequence of new flows $\{\widetilde M_{i,t}\}$ so that $\lambda_1(0)+\lambda_2(0)<\frac{1}{i}$ and $H(0,0)=1$ for all $i$.
$\medskip$

By the global convergence theorem (Theorem 1.12 of \cite{HK2}), after passing to a subsequence, the sequence of rescaled flows $\{\widetilde M_{i,t}\}$ converge locally smoothly to an $\alpha$-Andrews flow $\{M_{\infty,t}\}$ with convex time slices. And the limit flow satisfies $\lambda_1(0,0)=\lambda_2(0,0)=0$. By the strong maximum principle for tensors (see for example the appendix of \cite{W1}), the limit flow splits of a plane. By Fatou's lemma applied to each of the $F_{x_0, r}$ functionals individually, we see that the limit flow is also low entropy ($\lambda(M_{\infty, t}) < \Lambda_{n-2}$).
$\medskip$

Now take the blow-down of this limit flow at $t = -\infty$; by Huisken's monotonicity formula, we get a nontrivial (because $H(0,0)=1$)) self-shrinker which splits off a plane and which is mean convex. By the classification of mean-convex self-shrinkers \cite{CM} the entropy is then at least $\Lambda_{n-2}$; by the above argument though the blowdown should as well be low entropy, so we get a contradiction. 
\end{proof}

\subsection{Existence of small-entropy mean-convex cap}

In Haslhofer and Kleiner, when they perform surgery at necks they only have to ensure the resulting surface stays uniformly 2-convex and do not worry about the affect on entropy at all. But one can see by a straightforward computation (see the appendix) that in a toy model of the surgery similar to their construction, where a round cylinder $S^{n-1} \times \R$ is replaced with a half cylinder $S^{n-1} \times (-\infty, 0]$ and a cap, the entropy of the postsurgery model must be strictly greater than that of the round cylinder. 
$\medskip$

Estimating exactly how much the entropy increases directly seems to be nontrivial, even in the toy case. To eschew this problem, in this section we construct a cap model $C$ of low entropy, in a precise sense, by making use of the monotonicity of the entropy under the mean curvature flow.  To do this we construct a low entropy hypersurface, denoted below by $\Sigma$, that develops a neckpinch and is approximately cylindrical just away from the singularity in a precise way; as pointed out above, this example also shows that singularities are indeed a real possiblity for hypersurface $M^n$ with $\lambda(M) < \Lambda_{n-2}$. A timeslice right after the singularity then provides our cap. Namely, the main result in this section is the following:

\begin{prop}\label{constructioncap}
For any $\epsilon>0$, there exists $R_1$ such that if $R>R_1$, there exists a rotational symmetric n-dimensional cap model $C$, such that:
\begin{enumerate}
\item $C\subset B(0,4R)\subset\mathbb R^{n+1}$, 
\item $\lambda(C)\leq\Lambda_{n-1}+\epsilon$, 
\item$C\cap \R^{n+1} \setminus B(0,2R)$ agrees with a round half cylinder of radius 1 centered at the origin, and 
\item $C$ is mean-convex and $\alpha$ non-collapsed for some $\overline{\alpha} > 0$. 
\end{enumerate}
\end{prop}

Of course if $M$ is $\alpha$-noncollapsed for $\alpha > \overline{\alpha}$ it is also $\overline{\alpha}$-noncollapsed so the exact value of $\overline{\alpha}$ above is immaterial (although it will be close to that of a cylinder). Before proving the proposition we will need some lemmas. The first lemma says for some cases at least only $F_{x_0, r}$ of certain scales are relevent in the estimation of entropy. 

\begin{lem}\label{largeF}
For any surface $\Sigma\subset \mathbb R^{n+1}$ contained in $B^1(0,1)\times\mathbb \R^{n} \subset \R^{n+1}$ with $\lambda(\Sigma)\leq\Lambda_{n-2}$, there exists $r_1>0$ (depending on the growth rate and constant) such that
\begin{equation}
\begin{split}
\lambda(\Sigma)=\sup_{x_0\in \R^{n+1},r>0}F_{x_0,r}(\Sigma)=&\sup_{x_0\in B^1(0,1)\times \R^{n},r<r_1}F_{x_0,r}(\Sigma)\\
=&\sup_{x_0\in B^1(0,1)\times  \R^{n},r<r_1}\int_{\Sigma}\frac{1}{(4\pi r)^{\frac{n}{2}}}e^{\frac{-|x-x_0|^2}{4r}}d\mu_x\\
\end{split}
\end{equation}
namely the entropy will only be approximated by F-functionals on a bounded range of scales.
\end{lem}

\begin{proof}
That the supremum only needs to be taken with $x_0\in B^1(0,1)\times\mathbb R^{n}$ follows from lemma \ref{conhull} above and that the surface is supported in this solid round cylinder.
$\medskip$

By the entropy bound, we can get a uniform Euclidean volume bound on the surface $\Sigma$. $\mathrm{Vol}(\Sigma\cap B^{n+1}(p,r))\leq Cr^n$ for any $p,r$ and $C$ is a universal constant. The lemma follows if we can show $$\lim_{r\rightarrow\infty}\sup_{\R^{n+1}}\int_{\Sigma}\frac{1}{(4\pi r)^{\frac{n}{2}}}e^{\frac{-|x-x_0|^2}{4r}}d\mu_x=0$$

By breaking the integral up into integration on concentric annuli, it can be estimated by:
\begin{equation}
\begin{split}
&\int_{\Sigma}\frac{1}{(4\pi r)^{\frac{n}{2}}}e^{\frac{-|x-x_0|^2}{4r}}d\mu_x\\
=&\int_{\frac{\Sigma}{r}}\frac{1}{(4\pi  )^{\frac{n}{2}}}e^{\frac{-|x-x_0|^2}{4}}d\mu_x\\
=&\sum_{k=1}^\infty\int_{\frac{\Sigma}{r}\cap [B^{n+1}(0,k)\setminus B^{n+1}(0,k-1)]}\frac{1}{(4\pi )^{\frac{n}{2}}}e^{\frac{-|x-x_0|^2}{4}}d\mu_x\\
\leq&\sum_{k=0}^\infty  \left(\frac{C}{r} \right)^n\cdot e^{-(k-1)^2/4}\\
=&\frac{\widetilde C}{r^n}\\
\rightarrow&\text{ }0
\end{split}
\end{equation}
as $r\rightarrow\infty$. The volume bound $\mathrm{Vol}(\frac{\Sigma}{r}\cap [B^{n+1}(0,k)\setminus B^{n+1}(0,k-1)])\leq \frac{C}{r}$ is because after rescaling $\frac{\Sigma}{r}$ is contained in a round solid cylinder of radius $\frac{1}{r}$.
\end{proof}

In the next lemma we observe that the integral in the defintion of $F$-functionals is concentrated within a bounded set for a given bounded range of scales; essentially if the scales aren't let to be large the $F_{x_0, r}$ functionals must be concentrated near their basepoint $x_0$ by letting $C = \Lambda_{n-2}$ in the lemma below: 

\begin{lem}\label{farentropy}
For any $\epsilon, r_1, C_0>0$, there exists $R_0>>1$ such that if $R>R_0$, then for any $M^n\subset\mathbb R^{n+1}$ with entropy $\lambda(M)\leq  C_0 < \infty$ it's the case that: 
\begin{equation}
\sup_{x_0 \in \mathbb R^{n+1} ,r<r_1}F_{x_0,r}(M \cap B(x_0,R)^c)=\sup_{x_0\in \mathbb R^{n+1},r<r_1}\int_{M\setminus B^{n+1}(x_0,R)}\frac{1}{(4\pi r)^{\frac{n}{2}}}e^{\frac{-|x-x_0|^2}{4r}}d\mu\leq\epsilon
\end{equation}
\end{lem}
\begin{proof}
As above the entropy bound implies Euclidean volume bound and 
\begin{equation}
\begin{split}
&\sup_{r<r_1}\int_{\{M-x_0\}\setminus B^{n+1}(0,R)}\frac{1}{(4\pi r)^{\frac{n}{2}}}e^{\frac{-|x|^2}{4r}}d\mu\\
=&\sum_{k=1}^\infty\int_{\{M-x_0\}\cap (B^{n+1}(0,(k+1)R)\setminus B^{n+1}(0,k\cdot R)}\frac{1}{(4\pi r)^{\frac{n}{2}}}e^{\frac{-|x|^2}{4r_1}}d\mu\\
\leq&\sum_{k=1}^\infty C [((k+1)R)^n-(kR)^n]e^{-|kR|^2/(4r_1)}\\
\leq&\sum_{k=1}^\infty \widetilde C k^{n-1}R^ne^{-|k(R-1)|^2/4r_1}\cdot e^{-k^2(2R-1)/(4r_1)}\\
=&e^{-k^2(2R-1)/(4r_1)}\sum_{k=1}^\infty \widetilde C k^{n-1}R^ne^{-|k(R-1)|^2/4r_1}\\
\leq&\bar C e^{-k^2(2R-1)/(4r_1)}\\
\rightarrow&\text{ }0
\end{split}
\end{equation}
as $R\rightarrow\infty$.

So the lemma follows by choosing $R_0$ large enough.
\end{proof}

We also need to consider the following fact, which is a consequence of the continuity of each of the $F$ functionals having bounded gradient within a $C^3$ bounded family of submanifolds; see \cite{Mra1} section 5.
\begin{lem}\label{cpctsetentropy}
For any $\epsilon>0$, and $R>R_0$ chosen above, there exists $\delta(\epsilon,R)>0$ such that if $\overline C$ is the graph of $u$ over a round cylinder $C_{r_0}$ of radius $r_0$ centered at origin and  $||u||_{C_3(B(0,R))}\leq\delta$, then
\begin{equation}
|\lambda(\overline C\cap B(0,R))-\lambda(C_{r_0}\cap B(0,R))|<\epsilon
\end{equation}
\end{lem}

With this in mind we describe how to construct $\Sigma$. First, consider a part of a round cylinder ``threaded'' through a self shrinking torus. On either side, gradually start to let the cylinder flare out. Provided it doesn't change radius too quickly, by the lemma above its entropy will be very close to that of a cylinder. On the other hand, we can ensure very far away from the neck pinch that it will not be singular, because we will be able to fit very large spheres within it. Once we have flared the cylinders out enough to fit large spheres whose flows are still smooth before the shrinking torus flows to a point, let the cylinder radius level off (suitably gradually) and eventually cap it off on either end. The following schematic summarizes the construction: 

\begin{center}
$\includegraphics[scale = .75]{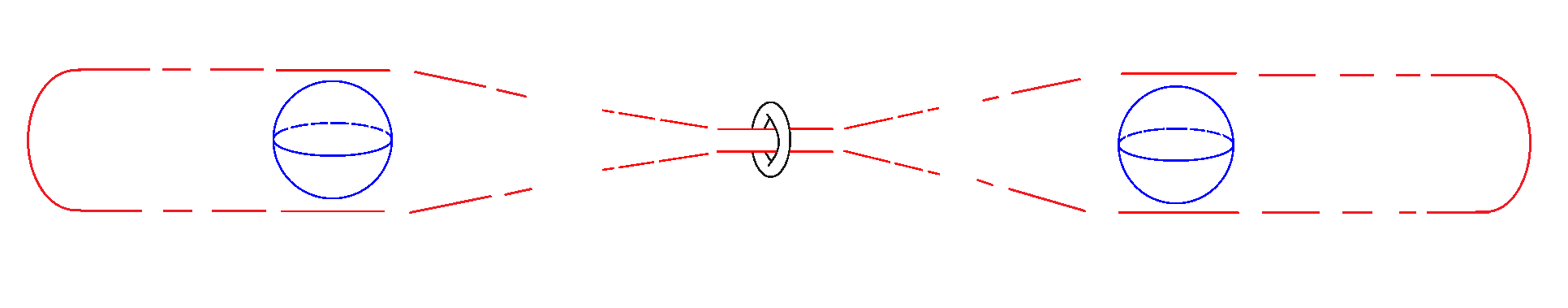}$
\end{center}

More precisely, let $\rho(x):\mathbb R\rightarrow[0,1]$ be a heavyside function, namely $\rho\in C_0^\infty$, $\rho(x)=0$ when $x\leq0$ and $\rho(x)=1$ when $x\geq1$. Let $m$ be chosen so that $||\frac{1}{m}\rho||_{C_3}\leq\delta$ in the condition of Lemma \ref{cpctsetentropy}. Let $W$ be the width of Angenent's shrinking torus at the time slice with inner radius 1. Define $\eta(x)=\frac{1}{m}\rho(x-2W)$. Denote the time $t_1 >0$ to be the time when the self shrinking torus of width $W$ shrinks to a point. 
$\medskip$

Define $\eta_k(x)=\frac{1+\sum_{j=1}^{km} (\eta(\frac{x}{2R_0}-j)+\eta(-\frac{x}{2R_0}-j))}{k+1}$ and choose $k$ large enough so that $k^2m^2>2nt_1$. Then the surface of revolution $\Sigma_r\subset\mathbb R^{n+1}$ defined by rotating the graph of $\eta_k$ around the $x_1$ axis must develop a neck-pinch singularity by the comparison principle for the mean curvature flow (as described in the background material, this is a consequece of the maximum principle). 
$\medskip$

This surface $\Sigma_r$ is contained in a solid round cylinder of radius 1 because $\eta_k(x)\leq1$ and it agrees with the round cylinder of radius 1 outside the ball of radius $2R_0(2W+mk+1)$.  So by Lemma \ref{largeF}, the entropy of $\Sigma_r$ are only approximated by $F$ functionals with bounded scales. Moreover,  by our choice of $R_0,m$, using lemma \ref{farentropy} and lemma \ref{cpctsetentropy}, we have $\lambda(\Sigma_r)\leq\Lambda_{n-1}+\epsilon$. 
$\medskip$

Now for any $\widetilde R>> R_1,R_0$, choose $R'>>10 \widetilde R$ and cap off $\Sigma_r$ by spherical caps outside the ball of radius $\widetilde R$ to get $\widetilde\Sigma_r$, which is of the shape of a long ``pill'' as in the diagram above. It will be strictly mean convex and hence $\alpha$ noncollapsed for some $\alpha$, and by the maximum principle this will be preserved under the flow. The $F_{x_0, r}$ functionals with $r < r_1$ and $x_0 \in B(0, \widetilde R)$ will be bounded by $\Lambda_{n-1} + \epsilon$ if $R'$ is sufficiently large since the bound holds for these functionals on $\Sigma_r$, which we will assume in the proceeding. 
$\medskip$

By the lemma below, which one can interpret as a pseudolocality result of sorts, if $R'$ is large enough, the evolution of $\widetilde\Sigma_r$ will be as close as we want to the evolution of a round cylinder in the annuli $A_{\widetilde R} = B(0,4\widetilde R) \setminus B(0, 2 \widetilde{R})$ which will let us control the geometry of the ``end'' of the cap: 
\begin{lem}\label{pl} Suppose $M_1$, $M_2$ are two submanifolds of $\R^N$ with entropy bounded by $\lambda$, whose mean curvature flow exists on the interval $[0,T]$ and $|A|^2$ is uniformly bounded initially by say $C$ in $B(0,R')$. Picking $\epsilon$ and $R$, there exists $R'(\epsilon, R, C, \lambda) > R$ so that if $M_1 \cap B(0,R') = M_2 \cap B(0, R')$ then $(M_1)_t \cap B(0,R)$ is $\epsilon$ close in $C^2$ local graphical norm to $(M_2)_t \cap B(0,R)$ for all $t \in [0,T]$. \end{lem} 
\begin{proof}
Without loss of generality $R= 1$. Suppose the statement isn't true; then there is a sequence of hypersurfaces $\{M_{1i}, M_{2i} \}$,  $R_i \to \infty$ and times $T_i \in [0,T]$ so that $M_{1i} = M_{2i}$ on $B(0,R_i)$ but $||M_{T_i} - {M_i}_{T_i}||_{C^2} > \epsilon$ in $B_0(1)$. By passing to subsequences by Arzela-Ascoli via the curvature bounds and area bounds we get limits $M_{1\infty}$, $M_{2\infty}$ so that $M_{1\infty} = M_{2\infty}$ (the flows of these manifolds will exist on $[0,T]$) but the flows don't agree at some time $T_1 \in [0,T]$; this is a contradiction by the uniqueness theorem for immersed MCF of Chen and Yin \cite{CY}. 
\end{proof}

Now we can give the construction of the low-entropy cap $C$:
\begin{proof}(of Proposition \ref{constructioncap})
A result by Angenent, Altschuler, and Giga \cite{AAG} on rotationally symmetric level set flows ensures that the singular times of the level set flow starting from $\widetilde\Sigma_r$ are discrete, so immediately after the neckpinch time $t_{neckpinch} < t_1$ (and because the entire surface doesn't go singular before $t_1$) our surface will be consist of two smooth components. In addition, the flow is nonfattening, and one can see that the smooth points will move by their mean curvature vector at all times. Our choice of cap model then is one connected component of a time slice immediately after the first neck-pinch singularity that lies inside the ball $B(0,4\widetilde R)$. By a result due also to Haslhofer and Kleiner (see theorem 1.5 in \cite{HK2}) the post-singular surface will be $\alpha$-noncollapsed as well. In parabolic neighborhoods not including singular points, the level set flow will be a mean curvature flow. It's clear by the rotational symmetry and mean convexity of the surface created that singularities only occur where a neckpinch occurs about the axis of symmetry and hence may be ruled out away from the origin for our time frames under consideration by a comparison argument using spheres, so the level set in a fixed open set not containing the origin will be a mean curvature flow in the classical sense with bounded curvature even through the first singular time. 
$\medskip$

Applying lemma \ref{pl} to balls covering $A_{\widetilde R}$ (if need be taking $\widetilde{R}$, and hence $R'$, larger), we see that the evolution of $\widetilde\Sigma_r$ is as close as we want in $C^2$ norm to the evolution of a round cylinder in $A_{\widetilde R}$ through the time of the neckpinch near the origin. By deforming it in $B(0,4\widetilde R)\setminus B(0,\frac{3}{2}\widetilde R)$, we can make it agree with a round cylinder in $B(0,4\widetilde R)\setminus B(0,2\widetilde R)$ and keep the entropy bound $\Lambda_{n-1}+\epsilon$ by lemma \ref{cpctsetentropy} above. Finally we extend this hypersurface by a half cylinder to get our cap $C$ -- again since entropy is only attained on bounded scales and the entropy of the cylinder is $\Lambda_{n-1}$ this will preserve the entropy bound of $C$ from the previous sentence. 
\end{proof}

To describe how we glue it in, note an upshot of the canonical neighborhood theoem above is that if the mean curvature at the locations we intend to do surgery is large enough, after rescaling to make the mean curvature one the surface will be as close as we want (in $C^3$ norm, say) in as large a neighborhood as we want to a round cylinder of radius one. Meanwhile, without loss of generality (by applying a suitable rescaling) our surgery cap candidate constructed in the previous subsection agrees with a round cylinder far enough away from the origin.
$\medskip$

The locations that we intend to do surgery at will have $H = H_{neck}$, as in \cite{HK2}, so choose $H_{neck} > H_{can}(\epsilon)$ with $\epsilon$ so that $1/\epsilon >> 2\widetilde{R}$ and $\epsilon < \delta/2$. Denote the rescaled flow about the surgery spot by $\widetilde{M}=\frac{M-p}{H_{neck}}$, by our choice of $\epsilon$ let us perform the cap gluing by smoothly transitioning from $\widetilde{M} \cap (B(0, 3\widetilde{R}) / B(0, 2\widetilde{R}))$ to $C \cap (B(0, 3\widetilde{R}) / B(0, 2\widetilde{R}))$. In particular the surgery only will change the hypersurface in the region $\widetilde{M} \cap (B(0, 3\widetilde{R})$ for the rescaled flow. Following the notation of Haslhofer and Kleiner we denote the surfaces pre and post-gluing by $\widetilde{M}^-$ and $\widetilde{M}^+$ respecitvely for the rescaled surface and $M^-$, $M^+$ for the original (spatial) scaled surfaces. 
$\medskip$

With regards to proposition 3.8, $\beta$ 2-convexity isn't strictly necessary (its only included in \cite{HK1} to preserve apriori curvature conditions) so we ignore that condition. Since the transition is taken where both surfaces are nearly cylindrical, $\alpha$-noncollapsing is preserved (although the apriori noncollapsing may need to be adjusted if the cap's is lower -- of course one may freely lower the constant $\alpha$ in the definition of noncollapsing). Items (1), (2) and (4) are clear as well; for the third point we note we may slightly bend the cylinder inwards without affecting the entropy in light of lemma 3.5 to make the postgluing domain satisfy (3).

\subsection{Estimation of entropy across surgeries.}

Now that we have the cap and how to glue it in, we analyze the change in entropy due to surgery and ensure that, if surgery parameters are picked correctly, the post surgery surface will still have low entropy. Morally speaking, since the cap was constructed to have entropy very close to that of the cylinder, the post gluing domain should have low entropy as well. The rub is that the contribution to the $F$ functionals near the surgery cap from the rest of the manifold could concievably be large, so that somehow even after placing caps the entropy is pushed over the low entropy threshold. We show with a careful choice of surgery parameters that this won't occur.
$\medskip$

There are two types of $F$-functionals to consider for us, those which are concentrated near $x_0$, or roughly when $r$ is small, and the diffuse ones where $r$ is roughly large. We start by showing we can find $c > 0$ so that all $F_{x_0, r}$ functionals with $r < c$ have $F_{x_0, r}(M_T) < \Pi$ (recall from the start of the section that we are considering just the single surgery time $T$). We then show by taking $H_{neck}$ large enough that we can arrange $F_{x_0, r}(M_T) < \Pi$ for $r > c$ as well. The choice of $c$ is somewhat delicate (at least in our approach) and we take an iterative approach to defining it to elucidate its choice and what other parameters it depends on: 
$\medskip$

First note as a byproduct of lemma \ref{farentropy}, which one can see by rescaling, for every $s, \epsilon_0 > 0$, there exists $c' = c'(\epsilon_0)$ so that if $r < c'$, $F_{x_0, r}(M) - F_{x_0, r}(M_T \cap B(x_0, s))  < \epsilon_0$, and we can apply it even after surgery caps are inserted because the entropy will certainly be at least finite and bounded by the area of $M_T$ (although our ultimate goal is a much sharper bound). Also by $\alpha$-noncollapsedness we know if $H$ has an upper bound $B$, $|A|^2$ does as well. Hence in sufficiently small balls it can made as close one wants to a plane (after rescaling), giving us as a consequence that for every $\epsilon_0, B > 0$ there is $s$ so that $F_{x_0, r}(M_T \cap B(x_0, s))  < 1+ \epsilon_0$ if $H < B$ in $B(x_0, s)$. Note trivially this conclusion is still true when lowering $c'$ but keeping $B$ fixed. 
$\medskip$

As a corollary of this observation we see that for every $\epsilon_1, Q > 0$ ($Q$ will be picked below), setting in the previous paragraph $B = 2QH_{can}(\epsilon_1)$, there is a $1 >> s_1 > 0$ so that if $F_{x_0, r}(M_T \cap B(x_0, s_1))  > 1+ \epsilon_0$, then $H(y) > 2QH_{can}(\epsilon_1)$ for some $y \in B(x_0, s_1) \cap M_T$.  Of course, the rough plan is to estimate the value of $F_{x_0, r}$ in terms of the canonical neighborhood models of these points in some manner, at least the ones that will be represented by parts of the surface that persist after surgery. 
$\medskip$

Considering $y$ (and $\epsilon_1$ and $Q$) as above then, there is an $s_2 << 1$ such that $H(x) \geq \frac{3Q}{2}H_{can}(\epsilon_1) > QH_{can}(\epsilon_1)$ for all $x \in B(y, s_2)$, for $s_2$ small enough by the gradient estimates (see theorem 1.10 in \cite{HK1}) applied at points where $H(y) = 2QH_{can}(\epsilon_1)$ -- if $H$ isn't strictly greater than this value at every point in $B(y, s_2)$ there must be a point where equality is met since $H$ is continuous. The important thing to note here is that $s_2$ only depends on $\epsilon_1$ (as $\epsilon_1$ decreases, so does $s_2$ because $H_{can}(\epsilon_1)$ increases) hence by taking $c_0 = c_0(s_2, \epsilon_0)$ (for which $r < c_0$) sufficiently small we may in fact arrange that if $F_{x_0, r}(M_T)  > 1+ \epsilon_0$ then $H > QH_{can}(\epsilon_1)$ in $M_T \cap B(x_0, s_2)$, where the ball considered is now centered at $x_0$, and so that $F_{x_0, r}(M_T) - F_{x_0, r}(M_T \cap B(x_0, s_2))  < \epsilon_0$.
$\medskip$

Intuitively speaking since we found one very high curvature point the point must be ``deep'' in the neck so only surrounded by high curvature points. With all this in mind set $c = c_0$, denote by $\overline{H}:= QH_{can}$, and as a first pass consider $F_{x_0, r}$ functionals such that the following hold; we will such $F$ functionals are very close to $\Lambda_{n-1}$:

\begin{enumerate}
\item $0 < r < c$ (note in the course of the proof, $c$ will be adjusted further).
\item $F_{x_0, r}(B(x_0, s_2) \cap M_T) > 1 + \epsilon_0$.
\item $B(x_0, s_2) \cap M_T$ only contains ``neck points'' in the sense discussed in section 2.2 above.
\item Furthermore, no surgeries are done in $B(x_0, s_2) \cap M_T$.
\end{enumerate}
Of course, if the second point is not satisfied then taking $\epsilon_0$ small enough $1 + 2\epsilon_0 < \Pi$, so these $F_{x_0, r}$ functionals will not potentially ruin the low entropy condition after surgery. We will discuss the complement of the other three cases below in the order of (4), (3), and then (1). 
$\medskip$

With regard to assumption (3), in fact we see points which flow to neck points by time $T$ and are already of high enough mean curvature to apply the canonical neighborhood theorem must themselves be neck points: if alternately they were modeled on caps, the flow would have already left $B(x_0, s_2)$ by time $T$ by taking $Q$ sufficiently large). With this in mind see by unpacking definitions the canonical neighborhood theorem for any point $x \in B(x_0, s_2) \cap M_t$ we have $H(x)\cdot(M-x)$ is $\epsilon_1$ close in $C^{\frac{1}{\epsilon_1}}$ topology to the round cylinder of radius $n-1$ in the ball $B(\frac{1}{\epsilon_1})$ provided $H > H_{can}(\epsilon_1)$.
$\medskip$

So, taking $s  = \min\{ \frac{10(n-1)}{\overline{H}}, s_2\}$ possibly less than $s_2$ and relabeling $s = s_2$, we have every point $M_t \cap B(x,s_2)$ is $1.1\epsilon_1$ close to a \textit{single} round cylinder (of radius 1) after rescaling (with our previous choice of $s_2$, instead this was only known to be true in every neighborhood of every point) when $H = \overline{H}$ somewhere within the ball (the coefficient 1.1 is to account that some of the points might not be at the center of the cylinder found, which is centered where $H = \overline{H}$, and the ball of radius $10(n-1)$ (after rescaling) used so whole cross sections of the cylinder near such points are contained in the ball). Note this means potentially lowering $c$ as well. 
$\medskip$

Now we estimate $F_{x_0,r}(M_T)$ (specifying pre or post surgery doesn't matter in this case by (4)), and we want to show under assumptions (1)--(4) the $F$ functionals must be small in the sense they can be bounded above by $\Lambda_{n-1}$ plus a small remainder. We first note the following lemma, a consequnece of the continunity and mean convexity of the flow and that without loss of generality for the initial surface $M$, $\sup H < .9\overline{H}$ for $\epsilon_1$ small enough (in other words $\gamma < .9\overline{H}$): 

\begin{lem}\label{cov}
Those points with mean curvature $H(x,t)>\overline{H}$ must be covered by union of balls $\bigcup\limits_{\tilde x, H(\tilde x,\tilde t)=\overline{H} \text{ for some $\tilde t\leq t$}}[B(\tilde x,s_2)]$.
\end{lem}


In the following few paragraphs we show, by adapting arguments of Gianniotis and Haslhofer \cite{GH} where they prove a bounded diameter theorem for 2-convex flows, that in these balls we find neck points whose axii are ``stable'' under the flow. The upshot is by time $T$ the surface will still be close to a (single) round cylinder in these balls, and hence have the associated $F$ functionals will have the bound as claimed: 
$\medskip$

By lemma \ref{cov} above, there are some points in the ball $B(x_0,s_2)$ of some previous time slice of the unscaled flow that with mean curvature exactly equal to $\overline{H}$ by our assumptions (1) - (4); denoting this time $t_1 < T$ we have by the argument in section 4.3 of \cite{GH} that every point in $M_{t} \cap B(x_0, s_2)$ is a $\epsilon_1$ strong neck at least when $Q \geq 10$ (see section 2 of their paper for this definition -- the point is this neck is closely modeled on a cylinder far back in time) for $t \in [T- t_1 - \tau, T-t_1]$ where $\tau = \frac{1}{16}H_{can}(\epsilon_1)^{-2}$, provided $\epsilon_1$ is small enough. By our choice of $s_2$ above, at time $t_1$ the model cylinder will be the same for all points in $B(x_0, s_2)$ (in other words, the axis of the model cylinder is fixed), relaxing $\epsilon_1$ to $1.1\epsilon_1$. 
$\medskip$

 Note by mean convexity and $\alpha$ noncollapsing that we'll continue to have $H > \widetilde{H}(\overline{H})$ along the space time track of points $p \in B(x_0, s_2)$ -- otherwise later on the flow will have a graphical piece ``bigger'' than the tube was at time $t_1$ contradicting set monotonicity of the flow. Its clear that as $\overline{H}$ increases $ \widetilde{H}$ does so by taking $Q$ larger we can arrange that $\widetilde{H}(\overline{H}) > 10H_{can}(\epsilon_1)$. Hence for all $T - t_1 \leq t \leq T$ we'll have every point in $M_t \cap B(x_0, s_2)$ is an $\epsilon_1$ strong neck point again by section 4.3 of \cite{GH}. Also note since this interval is closed we can cover it by finitely many subintervals of length $\tau$ and, by comparison with a sphere enclosing $M$, the number of subintervals is bounded uniformly irrespective of $T$ or $\overline{t}$ (since the flow will become null when the enclosing sphere does). 
$\medskip$

The point of all this is that we can then apply proposition 4.1 of \cite{GH}, which argues the axii of the model cylinders don't ``tilt'' much by using Colding and Minicozzi's Lojasewicz-Simon inequalities \cite{CM3}, in each of the intervals $[t - \tau, t]$ for $T - t_1 \leq t \leq T$ to conclude for $\epsilon_1' > 0$ that every point $p$ in the lemma above satisfies $M_T \cap B(p, s_2)$ is $\epsilon_1'$ close to a single cylinder (the one found at time $t_1$) if $\epsilon_1$ is sufficiently small. 
$\medskip$

Hence these $F$ functionals under consideration restricted to $M_{T} \cap B(x_0, s_2)$ are at most $\Lambda_{n-1}  + \epsilon_2$, where $\epsilon_2$ is a function of $\epsilon_1'$ (and hence $\epsilon_1$) which tends to zero as $\epsilon_1$ does. Thus:
\begin{equation}
\begin{split}
&\int_{M_{T}}\frac{1}{(4\pi r)^\frac{n}{2}}e^{\frac{-|x-x_0|^2}{4r}}\\
<& \Lambda_{n-1} + \epsilon_0 + \epsilon_2\\
\end{split}
\end{equation}
$\medskip$

If $\epsilon_1$ (and hence $\epsilon_2$) and $\epsilon_0$ are sufficiently small, $\Lambda_{n-1} + \epsilon_0 + \epsilon_2 < \Pi$; for future reference we note tracing through the argument above that equivalently the inequality holds if $\epsilon_0$ is sufficiently small and $H_{can}$ is sufficiently large, so really only depends on $\epsilon_1$ being sufficiently small.
$\medskip$

 Now suppose a surgery is done in $B(x_0, s_2)$ and let $F_{x_0, r}$ be an $F$ functional still satisfying properties (1) - (3) above. Note we still have $F_{x_0, r}(M_T) < \Lambda_{n-1} +\epsilon_0 + \epsilon_2$ by the work above before a surgery is done. 
$\medskip$

With that in mind, the surgery occurs in a small ball $B(q, r_0)$ about the center of the surgery region, where $r_0 = \frac{5\widetilde{R}}{H_{neck}}$, and by the design of the surgery caps $|F_{x_0, r}(M_T^+ \cap B(q, r_0)) - F_{q, r}(M_T^- \cap B(q, r_0))| < \epsilon$. In fact we see the $F_{x_0, r}$ functional will potentially increase the most under the surgery if it is situated right at the center $q$ of the surgery region, at least prior to the deletion of the high curvature region (by the symmetry of the Gaussian distribution); of course the subsequent deletion of the high curvature regions will only decrease each of the $F$ functionals and hence the entropy. Thus after surgery is done, $F_{x_0, r}(M_T^+) < \Lambda_{n-1}  + \epsilon_0 + \epsilon_2 + \epsilon$ and if $\epsilon$ is taken small enough this will still be less than $\Pi$. 
$\medskip$


We claim that in fact without loss of generality only the base points with $B(x_0, s_2)$ containing neck regions possibly with fresh surgery caps are the only ones one needs to consider. Suppose that the ancient flow one finds does lay in the region connecting $H_{th}$ to a point where $H = H_{trig}$. Well it must not lay in the low curvature region, since without loss of generality $\overline{H} > H_{th}$, so we see it must be discarded after all the surgeries at the surgery time (we are considering them one at a time) are complete. 
 \begin{center}
$\includegraphics[scale = .5]{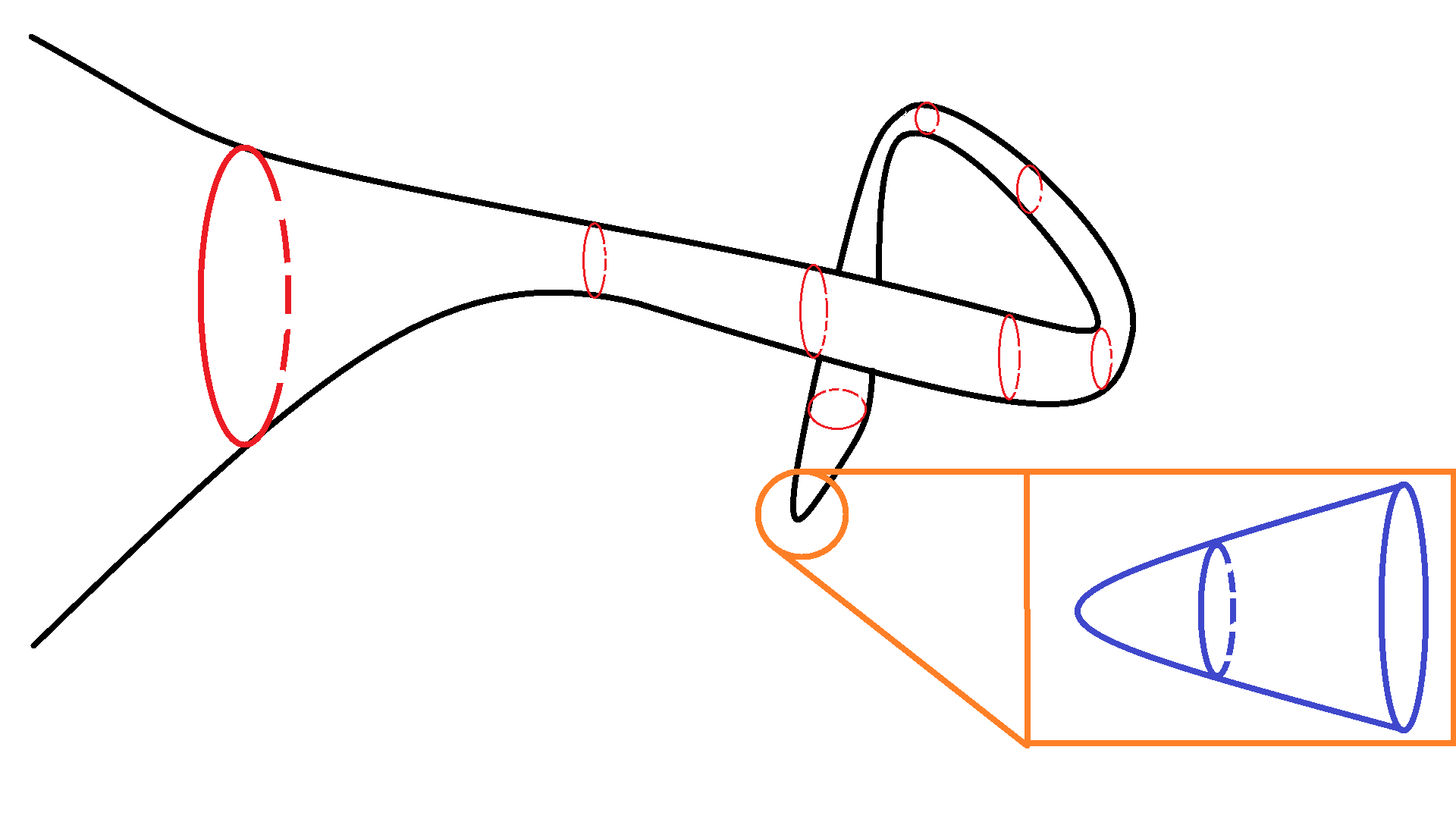}$
\end{center}
For example, in the diagram above (although in practice the circle might be quite a bit bigger relative to the scale of surgery), the circled tip of a high curature region is essentially modeled on a self translator and is not a neck point, but these are in the region of $M_T^\#$ that will be thrown away under surgery by the time the flow is allowed to continue again. 
$\medskip$

Note that the choice of $c$ ultimately only relies on $\alpha, \Pi,$ and $\gamma$ (initial curvature bound on the surface when start MCF with surgery), since we use the canonical neighborhood theorem. Concerning the surgery parameters $H_{th}, H_{neck}, H_{trig}$ we only required in the argument above that $\epsilon_1$ was small enough that $H_{th} < \overline{H}$, so since in the choice of surgery parameters one only requires $H_{trig}/H_{neck}, H_{neck}/H_{th}, H_{neck}$ be sufficiently large we may take $H_{neck}$ large freely in the next part of the argument, where we consider the case when $r \geq c$.  
$\medskip$ 

In fact, in this case we will see the $F$ functionals are in fact nonincreasing if $H_{neck}$ is taken large enough. For organizational conveinence we consider two domains, one about the surgery region centered at $q$ and the other ``far'' from the surgery, which we denote $U_e$ and $U_f$ respectively. More precisely, let $U_e = B(q, \frac{5\widetilde{R}}{H_{neck}})$, $\widetilde{R}$ as specified above, and let $U_f$ be its complement. We see the surgery takes place entirely within $U_e$. 
$\medskip$

First we will show if the surgery region $U_e$ is sufficiently small the Gaussian weights of $F$ functionals with a lower bound on $r$ are nearly constant within it in a sense made precise below. Then to conclude we use the following observation corresponding to the cap having less volume than the cylinder:
\begin{lem}\label{vol}
For $\widetilde{R}$ defined previously in the construction of the cap, there is an $0 < \eta(\widetilde{R}) < 1$ depending only on the surgery cap and $R$ such that
$$Vol(\widetilde M_{T}^+ \cap B(0, 3\widetilde{R}) )< \eta Vol(\widetilde M_{T}^- \cap  B(0, 3\widetilde{R}) )$$
where $\widetilde M_{T}=\frac{M-q}{H_{neck}}$
\end{lem}
Where above $\widetilde U_e = H_{neck}(U_e - q)$ , $\widetilde U_f = H_{neck}(U_f - q)$, and $\widetilde M_T = H_{neck}(M_T - q)$ denote the rescaled versions of $U_e$, $U_f$, and $M_T$ where the surgery center $q$ is translated to the origin. 
$\medskip$

To begin, note that the gradient of the Gaussian weight $e^{\frac{-|x -x_0|^2}{4r}}$ of a $F_{x_0, r}$ functional is given by: 
 \begin{equation}
\nabla e^\frac{-|x - x_0|^2}{4r} =\frac{ -2(x-x_0)}{4r}  e^\frac{-|x - x_0|^2}{4r} 
\end{equation}
Since by lemma \ref{conhull} the entropy for a compact hypersurface will be attained by an $F$ functional centered in its convex hull, without loss of generality $x_0$ is in the convex hull of $M^-$.  Since for such $x_0$ we have $|x-x_0|\leq D < \infty$\footnote{of course, the diameter is decreasing under the flow so is uniformly bounded by the diameter of the initial time slice}, we see for a lower bound $c$ on $r$ we have $|\nabla e^\frac{-|x - x_0|^2}{4r}| \leq \frac{D}{2c} < \infty$ for any choice of $x \in M_T$. Denote this upper bound by $\rho$. 
$\medskip$

We also note similarly for $r > c$ that the Gaussian weight of a $F_{x_0, r}$ functional (with $x_0$ in the convex hull of $M$) is bounded below by $e^{\frac{-D^2}{4c}} > 0$; denote this lower bound by $\sigma$.  Also denote by $m_{x_0, r}$ and $M_{x_0, r}$ the minimum and maximum respectively of the Gaussian weight of $F_{x_0, r}$ in $U_e$. Then the following is true:
\begin{equation}\label{weightgrad}
1 \geq \frac{m_{x_0, r}}{M_{x_0, r}}\geq \frac{m_{x_0, r}}{m_{x_0, r}+r_e\rho} \geq \frac{\sigma}{\sigma  + r_e \rho} = 1 - \frac{r_e \rho}{\sigma + r_e \rho} 
\end{equation}
Since $\sigma > 0$ and $\rho <\infty$ we can make this quotient as close to one as we like by making $r_e$ sufficiently small; in other words we can make the ratio of the minimum to the maximum of the weight in these $F$ functionals as close to 1 as we want in $U_e$ by increasing $H_{neck}$. Switching to the translated and rescaled picture (the ratio persists under rescaling), we have for $x_0\in\widetilde U_f$ and for $r>c$ the following:
\begin{equation}
\begin{split}
&F_{x_0,r}(\widetilde M^+)\\
=&\int_{\widetilde M^+}\frac{1}{(4\pi r)^{\frac{n}{2}}}e^{\frac{-|x-x_0|^2}{4r}}\\
\leq&\int_{\widetilde M^+\cap \widetilde U_f}\frac{1}{(4\pi r)^{\frac{n}{2}}}e^{\frac{-|x-x_0|^2}{4r}}+\int_{\widetilde M^+\cap \widetilde U_e}\frac{1}{(4\pi r)^{\frac{n}{2}}}e^{\frac{-|x-x_0|^2}{4r}}\\
&\text{(because surgery only happens in $\widetilde U_e$)}\\
=&\int_{\widetilde M^- \cap \widetilde U_f}\frac{1}{(4\pi r)^{\frac{n}{2}}}e^{\frac{-|x-x_0|^2}{4r}}+\int_{\widetilde M^+\cap \widetilde U_e}\frac{1}{(4\pi r)^{\frac{n}{2}}}e^{\frac{-|x-x_0|^2}{4r}}\\
\leq&\int_{\widetilde M^- \cap \widetilde U_f}\frac{1}{(4\pi r)^{\frac{n}{2}}}e^{\frac{-|x-x_0|^2}{4r}}+\int_{\widetilde M^+\cap \widetilde U_e}\frac{1}{(4\pi r)^{\frac{n}{2}}} M_{x_0, r}
\end{split}
\end{equation}
\begin{equation}
\begin{split}
&\text{(by Lemma \ref{vol})}\\
\leq&\int_{\widetilde M^- \cap \widetilde U_f}\frac{1}{(4\pi r)^{\frac{n}{2}}}e^{\frac{-|x-x_0|^2}{4r}}+\int_{\widetilde M^- \cap \widetilde U_e}\frac{1}{(4\pi r)^{\frac{n}{2}}}\eta M_{x_0, r}\\
&\text{(from the discussion after \ref{weightgrad} and taking } H_{neck} \text{ large enough})\\
\leq&\int_{\widetilde M^-\cap \widetilde U_f}\frac{1}{(4\pi r)^{\frac{n}{2}}}e^{\frac{-|x-x_0|^2}{4r}}+\int_{\widetilde M^-\cap \widetilde U_e}\frac{1}{(4\pi r)^{\frac{n}{2}}}e^{\frac{-|x-x_0|^2}{4r}}\\
=&F_{x_0,r}(\widetilde M^-)\\
\end{split}
\end{equation}
$\medskip$
So that these $F$ functionals don't increase under the surgery as claimed.  In all cases then we see the $F$ functionals $F_{x_0, r}$ either didn't increase after the surgery or they are bounded after the surgery by $\Lambda_{n-1} + \epsilon_0 + \epsilon_3 + \epsilon$ from above, with a prudent choice of surgery parameter $H_{neck}$ and sufficient choice of $\epsilon_i$. Picking $\epsilon_i$ sufficiently small gives $\Lambda_{n-1} + \epsilon_0 + \epsilon_2 + \epsilon < \Pi < \Lambda_{n-2}$ (we stipulated $\Lambda_{n-1} < \Pi$). To reiterate crucially we see from above we fixed these parameters and $\textit{then}$ varied $H_{neck}$ appropriately, without ruining our choice of $\epsilon_i$: this is because we only required along the way that $H_{th} < \overline{H}$, and we need to be able to take $H_{neck}$ large and luckily not small. In total the postsurgery surface can be arranged to be low entropy and we are done.

\section{Application to Self Shrinkers: Proof of Theorem 1.4}
In this section we show how our constructed mean curvature flow with surgery for mean convex hypersurfaces can be used to study self shrinkers of low entropy by considering a different (from $H > 0$) convexity condition that is preserved under the flow -- this condition is morally mean convexity for the renormalized mean curvature flow (abbreviated RMCF), which is related to the ``regular'' mean curvature flow by a reparameterization. In this section for the most part we do not work directly in the RMCF though (although it is used in an essential way in a step). Our starting point is the following observation:
\begin{lem} Closed self shrinkers that aren't already round may be perturbed to be $2H - \langle x, \nu \rangle$ $\alpha$ non collapsed, for some $\alpha > 0$, in an entropy nonincreasing way. 
\end{lem} 
Here naturally we say a hypersurface is $2H - \langle x, \nu \rangle$ $\alpha$-noncollapsed if it is $\alpha$-noncollapsed in the sense of section 2.1 above, except with respect to the quantity $2H -  \langle x, \nu \rangle$ instead of $H$. This observation is essentially lemma 1.2 in \cite{CIMW},  where one perturbs by the first eigenfunction of the stability operator $L$ introduced above - the noncollapsedness part then follows from compactness of the perturbed self shrinker. Note also that doing such a slight perturbation doesn't change topology and of course if the self shrinker is already round and compact it must be a sphere, so we have nothing to do. We intend to run the flow with surgery on these perturbations. 
$\medskip$

The perturbed surface starting from time $s=-1$ (as is traditional with self shrinkers) will become extinct before $s=0$ by comparison with the original shrinker. To stay in line with the rest of the article we shift the time $t=s+1$ so the flow of the perturbed surface will exists in $t\in[0,1)$ though. Following Lin \cite{L} (see also \cite{Smk}), we see that the quantity $F =  (2-2t)H - \langle x, \nu \rangle$ satisfies
\begin{equation}
\frac{dF}{dt} = \Delta F + |A|^2 F
\end{equation}
and hence $F$ $\alpha$-noncollapsing is preserved under the flow; when $t = 0$ this is exactly that the surface is $2H - \langle x, \nu \rangle$ $\alpha$-noncollapsed.  Since this family is our main interest of study in this section we formally define it:
\begin{defn} Denote by $\Sigma = \Sigma(\alpha, C, D, \Lambda)$ as the set of hypersurfaces:
\begin{enumerate}
\item $M \in \Sigma$ is initially $2H - \langle x, \nu \rangle$ $\alpha$-noncollapsed
\item $|A|^2 < C$ on $M$
\item $Diam(M) \leq D$
\item The level set flow starting from $M$ is empty strictly before $t = 1$. 
\item $\lambda(M) < \Pi$, where $\Lambda_{n-1} < \Pi < \Lambda_{n-2}$. 
\end{enumerate}
\end{defn}
We see any perturbed self shrinker above will be in $\Sigma$ for some choice of parameters. Our goal of this section rephrased then is to show existence of the mean curvature flow with surgery out of elements of $\Sigma$ for any choice of $\alpha >0$, $C > 0$, $\sigma < 1$, $D <\infty$ and $\Lambda_{n-1} <\Pi < \Lambda_{n-2}$. 
$\medskip$

We will say a point $p$ is $F$ $\alpha$-noncollapsed if it admits inner and outer osculating spheres of radius $\frac{\alpha}{F(p)}$. To show the existence of a surgery flow, we will show points of high curvature are $H$ noncollapsed and appeal to the mean curvature flow with surgery as already defined above for mean convex, low entropy mean curvature flow\footnote{As a (non rigourous) motivation, morally if a point has high curvature $H$ should be large, so since $t < 1$ and if $(2 - 2t)H - \langle x, \nu \rangle > 0$ then $H$ should be positive. The ``moral'' is true in our case due to the $F$ noncollapsing.}. Of course some details need to be checked;  to start one would want a uniform lower bound on $H$ for which we know the surface will be $H$ $\alpha$-noncollapsed, and one would only want the $H$ noncollapsing constant, which we'll denote $\hat\alpha$, to only depend on the parameters above describing the set $\Sigma$. This brings us to our first lemma:
\begin{lem} 
Suppose $M \in \Sigma$, and suppose $p \in M_t$ has $|A|^2(p) > n\frac{9D^2}{ \alpha^2}$. Then $H(p) > D$ and $p$ is $H$ $\hat\alpha$-noncollapsed for $\hat\alpha = \frac{\alpha}{3}$
\end{lem}
\begin{pf} 
To see this, we note by the $(2-2t)H - \langle x, \nu \rangle$ $\alpha$-noncollapsing that, denoting by $\lambda_i(p)$ the $i$-th principal curvature of $A$ at $p$:
\begin{equation}
|\lambda_i| \leq \frac{(2-2t)H - \langle x, \nu \rangle}{\alpha}
\end{equation} 
We then trivially estimate the numerator using the diameter of $M$ is initially bounded by $D$ and this persists under the flow:
\begin{equation}
(2-2t)H - \langle x, \nu \rangle \leq (2 -2t)H + D
\end{equation} 
Recalling that $|A|^2(p)$ is the sum of the squares of the principal curvatures of $M_t$ at $p$, if $|A|^2 >n \frac{9D^2}{ \alpha^2}$ then there must be some principal curvature $\lambda_j(p)$ so that $\lambda_j^2(p) > \frac{9D^2}{\alpha^2}$. Putting (4.2) and (4.3) together then yields:
\begin{equation} 
2D\leq |\lambda_j|\alpha - D\leq (2-2t)H
\end{equation} 
Since $0 \leq t < 1$ then $H > D$. Hence we get mean convexity; to get the statement on osculating spheres note at such points $p$:
\begin{equation} 
(2-2t)H - \langle x, \nu \rangle \leq 2H + D \leq 3H
\end{equation} 

From this because there are inner and outer osculating spheres at $p \in M_t$ of radius $\frac{\alpha}{(2-2t)H - \langle x, \nu \rangle}$, there are inner and outer osculating spheres of radius $\frac{\alpha}{3H}$ completing the proof. 
\end{pf}
As a corollary of this we obtain the following:
\begin{lem} 
Suppose $M \in \Sigma$. Then if $H > n^{3/2} \frac{3D}{ \alpha} =: \Phi$, $p$ is $H$ $\hat\alpha$-noncollapsed. 
\end{lem}
\begin{pf}
Suppose $H >  n^{3/2}\frac{3D}{ \alpha}$. Then one of the principal curvatures $\lambda_j > \sqrt{n} \frac{3D}{\alpha}$, which then implies $|A|^2 > n\frac{9D^2}{ \alpha^2}$.
\end{pf}

Thus points where $H$ is sufficiently large will be noncollapsed in the typical sense (i.e. $H$ noncollapsed). It is clear from the local nature of the proofs in \cite{HK1} that the canonical neighborhood theorem will still be true at points $p$ with some uniformly sized parabolic ball $P(p,t, \sigma)$ about them and will yield an $H_{can}(\epsilon)$ with $\frac{1}{H_{can}(\epsilon)} < \sigma$ (for a given choice of $\epsilon > 0$), the intuitive reason being that it is a blowup argument and the points outside the parabolic neighborhood will be rescaled to spacetime infinity. 
$\medskip$

Note that the following is not a ``surgery version,'' i.e. is only stated up to the first singular time; we will explain below what changes are necessary after the first surgery time, and $\Phi$ is as in the lemma above: 
\begin{prop}\label{MCnbhd} Let $T_0$ be some time before the first singular time and suppose $M \in \Sigma$ and $p \in M_t$ is so that $H(p)  = 2\Phi$ and $t > T_0$. Then there exists $\sigma > 0$ so that in $P(p,t, \sigma)$, $\Phi < H < 3 \Phi$ and so all points $q \in P(p,t, \sigma)$ are $\alpha$-noncollapsed. Furthermore if $(p,t)$ is such that $H > 2 \Phi$ and $t > 2T_0$, then $H > \Phi$ for all $q \in P(p,t, \sigma)$. 
\end{prop}
\begin{pf}Before starting we note that the time until the first singular time is uniformly bounded below by the evolution equation for $|A|^2$ and the uniform initial curvature bound $C$. First suppose that $H(p) = 2\Phi$; from the $\ell = 1$ curvature estimates (this is where having a nonempty parabolic ball coming from the stipulation $t > T_0$ is used) we immediately obtain a ball (i.e. for the fixed time slice $t$) in which within $B(p, \mu)$ $3\Phi/2 < H < 5\Phi/2$. Now define $s$ be the infimum of $|t - t'|$ over the set of times $t'$ before $t$ in which the spacetime neighbrohood $B(p, \mu) \times [t', t]$ contains a point where $H(p) = \Phi$; we will clearly be done if we can show $s > 0$. Recalling the following basic evolution equation: 
\begin{equation}
\frac{dH}{dt} = \Delta H + |A|^2 H
\end{equation}
By the $\ell = 2$ and $\ell = 0$ local curvature estimates, this is clearly bounded uniformly (in terms of $D$ and $\alpha$). Thus by integrating we see $s > 0$. 
$\medskip$

Now suppose that $H(p) > 2\Phi$. Then by the continuity of $H$ there is a backwards spacetime neighborhood $U$ of $p$ in which $H(q) = 2\Phi$ on $\partial U$ and $H(p) > 2\Phi$ in the interior of $U$. Consider $\overline{U} = U \cap \{t > T_0 \}$ (this could very well be just $U$); if $\partial \overline{U}$ consists of only points $q$ with $H(q) = 2\Phi$ we get the result from the above; the other case is when there are boundary points with $H(q) = 2 \Phi$. But since the claim is for points $(p,t)$ with $t > 2T_0$ , and from the proof above $\sigma^2 < T_0$, such points are sufficiently far away in the past that the statement holds. 
\end{pf}
$\medskip$

Before moving on to describing the construction of the surgery flow we briefly discuss the aforementioned local curvature estimates for $F$-noncollapsing flows- these were used in the above proposition. Lin showed these for starshaped mean curvature flow (theorem 3.1 in \cite{Lin}), although he remarks (specifically see remark 2.4 in \cite{Lin}) that these estimates are true for a wide class of flows including ours: 
\begin{thm}\label{localcurv} (Local curvature estimate). There exist $\rho = \rho(\alpha, \beta) > 0$ and $C_\ell = C_\ell(\alpha, \beta) < \infty$ so that if $M_t$ is a mean curvature flow with initial condition in $\Sigma$ defined in a parabolic ball $P(p,t,r)$ with $H(p,t) \leq r^{-1}$, then $M_t$ is smooth in the parabolic ball $P(p,t, \rho r)$ and 
\begin{equation}
\sup\limits_{P(p, t, \rho r)} | \nabla^\ell A| \leq C_\ell r^{-(\ell + 1)}
\end{equation}
\end{thm} 
Above $\beta$ is a lower bound on $H$ which for our case follows easily from the diameter bounds on the initial surface. The proof of this goes exactly as in the proof of theorem 3.1 in \cite{Lin}, where Lin proves it for starshaped mean curvature flow. The only point of that proof that might require some clarification is how to check his claim 3.8 (a one-sided minimization result that allows one to upgrade Hausdorff convergnece of a certain sequence in the proof to smooth convergence on compact sets).
$\medskip$
 
Here we briefly reparameterize to the RMCF.  To begin note that by setting $s =-(1-t)$, then the initial data is the $t = -1$ time slice of a flow $M_t$ defined on the time interval $[-1, -(1 - \sigma))$ and, in this parameterization, we have 
\begin{equation}
-2s H  - \langle x, \nu \rangle > 0
\end{equation}
We rescale the flow as follows:
\begin{equation}
\widetilde x(\cdot,\tau)=\frac{1}{\sqrt{-s}}x(\cdot,s), \text{ } \tau=-\log(-s)
\end{equation}
where $s\in[-1,0), \tau\in[0,+\infty)$.

The mean curvature is rescaled by $\widetilde H=\sqrt{-s}H$, and the rescaled flow satisfies the rescaled mean curvature flow equation (see \cite{H} for detail)
\begin{equation}
\begin{split}
(\frac{\partial}{\partial\tau}\widetilde X)^\perp &= -\widetilde H\widetilde \nu + \frac{\langle \widetilde x, \widetilde\nu  \rangle}{2})\widetilde\nu\\
&=(-\sqrt{-s}H+\frac{1}{2\sqrt{-s}}\langle x,\nu \rangle)\nu
\end{split}
\end{equation}
one easily sees the speed is negative: 
\begin{equation}
\begin{split}
&(-\sqrt{-s}H+\frac{1}{2\sqrt{-s}}\langle x,\nu \rangle)\\
=&(-\frac{1}{2\sqrt{-s}})(-2s\cdot H-\langle x,\nu \rangle)\\
=&(-\frac{1}{2\sqrt{-s}})F\\
<&0\\
\end{split}
\end{equation}
With this in mind we then define (again slightly different from Lin) the natural weighted area, where $N$ is a hypersurface:
\begin{equation}
Area_w(N) = \int_N e^{\frac{-|x|^2}{4}} d\mu
\end{equation} 
One can easily then check using this weighted area that the proof of claim 3.8 in \cite{Lin} goes through. Now we are ready to describe how the flow with surgery is constructed. From the curvature bounds (2) in definition 4.1, there is a uniform lower bound $\overline{T}$ for which the smooth flow exists. Let $T_0$ be so that $2T_0 < \overline{T}$ We thus obtain a $\sigma$ as in propostion \ref{MCnbhd}. 
$\medskip$

So take $M \in \Sigma$. From theorem 1.1 for our $\hat\alpha = \alpha/3$, using the small modification of the canonical neighborhood theorem as described before proposition \ref{MCnbhd}, we obtain $H_{th}, H_{neck}, H_{trig}$ for, it were true $M \in \mathcal{M}(\hat\alpha, n,\Lambda)$ (which in our case it isn't, of course), there would be a flow with surgery starting at $M$. We will then argue in fact though that one can define a flow with surgery for these choices of parameters for an $M \in \Sigma$ if the following conditions hold. The first condition can always be arranged because we have the freedom to take $H_{th}$ as large as we like. If the second condition is not true for any $H_{trig} > 0$ then lemma 4.2 gives the first singular time is less than $2T_0$, which contradicts that $2T_0 < \overline{T}$. On the other hand by $F$ noncollapsing and that $(2 -2t) > 0$ (by (4) in definition 4.1) we must have $H \to \infty$ as $t$ approaches the first singular time, so we see the second condition holds if $H_{th} < H_{neck}$ is large enough. The third condition is to ensure that the canonical neighborhood theorem can be used as described before the proof of proposition \ref{MCnbhd} and, again, can be arranged by potentially taking $H_{th}$ larger.
\begin{enumerate}
\item $H_{th} > 3 \Phi$
\item the first time $T_1$ for which $H = H_{trig}$ somewhere is greater than $2T_0$.
\item $\frac{1}{H_{th}} < \sigma$
\end{enumerate}
With this in mind start flowing $M$ by the mean curvature flow. Note if at any time along the flow no point satisfies $H = H_{trig}$, we will be able to continue the flow because we see from equation (4.2) that $|A|$ will be bounded in terms of $H$ (note this is a weaker statement that the $\ell = 0$ curvature estimates).
$\medskip$

Let $T_1 > 0$ be the first time that $H = H_{trig}$ somewhere on $M_{T_1}$. Then by proposition \ref{MCnbhd} we have in $P(p, T_1, \sigma)$ that the flow of $M$ is $H$ $\hat\alpha$-noncollapsed for all $q \in U_{H_{th}}(p) = \{x \in M_{T_1} \mid H(x) > H_{th} \}$, using of course assumption (1) on $H_{th}$ above. By items (2) and (3) the canonical neighborhood theorem holds in the region where $H > H_{th}$ by the time $T_1$ (for the $\epsilon$ necessary to do the surgery) so we can find ``separating points'' where $H = H_{neck}$ as usual; we cut the necks and place caps. The high curvature regions are of controlled topology and discarded. The mean curvature of the postsurgery domain is bounded by approximately $H_{neck}$, implying again from the $F$-noncollapsing a bound on $|A|$. Hence the flow can be restarted and ran smoothly a definite amount of time before $H = H_{trig}$ again. 
$\medskip$

There are still a couple things to check to have enough control over the high curvature region the next time $H = H_{neck}$ somewhere to do surgery again. First one needs to check that one can glue in the caps so as to preserve $F$ $\alpha$-noncollapsing and low entropy. We recall from the previous section that the cap is $\alpha$ noncollapsed for some $\overline{\alpha}$. Since the cap is placed at a point where $H > 2 \Phi$, namely as we see from lemma \ref{MCnbhd} so that it is larger than $2D$ ($\alpha < 1$, without loss of generality). Hence $2H - D > \frac{H}{2}$, so that the cap will also be $F$ $\alpha$-noncollapsed for $\alpha = \frac{\overline{\alpha}}{2}$. After potentially lowering, without loss of generality, the $F$ noncollapsing constant so that it is less than $\frac{\overline{\alpha}}{2}$ the surgery will thus be done so that it preserves $F$ noncollapsedness without degeneration of the constant. 
$\medskip$

In order to be able to do the surgery at the second time $T_2$ when $H = H_{trig}$ occurs we also need to check the validity of proposition \ref{MCnbhd} after $T_1$ - we see this boils down to showing the local curvature estimates when the relevent parabolic balls $P(p, t, \sigma)$, where $H(p) = 2 \Phi$, possibly contain surgeries from the first surgery time $T_1$ \footnote{for example, it is concievable the cap placed after a surgery is quickly ``pushed down,'' so that the curvature at the tip is low almost immediately after the surgery.}. For the $H$ noncollapsed surgery flows, this is essentially the content of theorem 1.6 in \cite{HK1} by Haslhofer and Kleiner, although their statement is significantly more general than we need because all surgeries for us occur with the same scale $s = \frac{n-1}{H_{neck}}$. 
$\medskip$

The proof in the $H$ noncollapsed case due to Haslhofer and Kleiner proceeds schematically the same as the smooth/no surgery version of the estimates, except there are now three cases dealing with whether or not a certain sequence of flows with surgery in the argument (it is a compactness-contradiction argument) have any surgeries and, if so, what types there are. Case (3) in their argument, if there are ``microscopic'' surgeries (see page 12 of \cite{HK1}) is far and away the most complicated part of the argument. Luckily this doesn't apply to us, because we see (since we only have to plan for surgeries happening at the scale $s = \frac{n-1}{H_{neck}} > 0$) we only have to deal with cases (1) and (2). The arguments there are very quick and clearly don't use mean convexity; only pseudolocality and the local gradient estimates in the presence of no surgeries (which we already have) are needed. 
$\medskip$

Thus we continue the flow with surgery until the surface is exhausted -  which occurs in finite time (in fact, before $t = 1$). We obtain the following topological consequence of the flow with surgery. 
\begin{cor} Let $M^n$ be a compact self shrinker with entropy less that $\Lambda_{n-2}$. Then $M$ is diffeomorphic to either $S^n$ or a connect sum of $S^{n-1} \times S^1$. 
\end{cor} 
At this step we compare with the result of Hershkovits and White: 
\begin{thm} (Theorem 1.1 in \cite{HW}) Suppose that $M \subset \R^{n+1}$ is a codimension-one, smooth, closed self shrinker with nontrivial $k^{th}$ homology. Then the entropy of $M$ is greater than or equal to the entropy of a round $k$-sphere. If equality holds, then $M$ is a round $k$-sphere in $\R^{k+1}$. 
\end{thm}

From the corollary, we see that all of the possible (apriori) topological types of compact self shrinkers except those diffeomorphic to $S^n$ have nontrivial first homology class. But for $n > 3$, $\Lambda_{n-2}$ is strictly less than $\Lambda_1$, so the only possible topological type of low entropy self shrinker for $n > 3$ is $S^n$. When $n =3$ the rigidity condition says the flow must be round $S^1 \subset \R^2$, which is another contradiction. Hence all such surfaces are diffeomoprhic to $S^n$ (with the standard differentiable structure). 
$\medskip$

To gain the isotopy one wishes to argue as in the inductive proof of theorem 1.3 in \cite{Mra} by ``extending the necks" found by the surgery (the induction is on number of surgeries encountered). In the mean convex case what one does is consider the future paths, which we'll refer to below as strings, traced out by the tips of the surgery caps under the flow with surgery and ``thicken'' (take a small tubular neighborhood) the paths to gain an isotopy to a tubular neighborhood of an embedded tree (i.e. a graph with no cycles) which one may then contract down to the round sphere. 
$\medskip$

Things are a little more complicated in our case however because the flow is not monotone since we don't have (traditional) mean convexity. Refering to the notation of proposition 3.2 in \cite{Mra}, one would hope to be able to intermitently stop the isotopies of $\mathcal{A}$ and $\mathcal{B}$ to isotope the strings out of the way. There is a topological obstruction to this however, which we see comes from the simple connectedness (or lack thereof) of the complement of the region bounded by the flow\footnote{For example, although this situation couldn't happen in our case, if the string was threaded through a low curvature region left after the surgery that eventually shrunk to a point it couldn't be moved out of the way.}
$\medskip$

For us though we know already our original hypersurface is diffeomorphic to $S^n$ (and hence any pieces leftover after the flow). By the noncollapsing, the low curvature pieces will be locally flat (in that there is a small tubular neighborhood of them will also be embedded), and using this one could appeal to the generalized Schoenflies theorem of Mazur \cite{Maz} and Brown \cite{Bro} to get an isotopy through locally flat embeddings. 
$\medskip$

One can argue alternately in a simpler way as follows though, with the added benefit of the isotopy being through smooth embeddings. By reparameterizing the flow to the renormalized mean curvature flow as alluded to above (see \cite{Mra2} for a more in depth discussion of surgery in the context of the renormalized mean curvature flow), the flow with surgery above is shrinker mean convex and hence ``monotone.'' This way there will be no issue with the strings (constructed after applying the reparameterization) ``bumping'' into the future flows of low curvature regions, and so we gain an isotopy by just following the argument of theorem 1.3 in \cite{Mra} directly.

\section{Concluding Remarks.}

There are other convexity assumptions under which a surgery theory are almost certainly possible; as noted in \cite{Lin} $F$ $\alpha$-noncollapsing is preserved for any $F$ of the form $F = a_1 \langle x, \nu \rangle + (a_2  + 2a_1t)H$ and the necessary estimates for surgery should go through (the case above is $a_1 = -1, a_2 = 2$). Note also that when $n = 2$, the flow with surgery is possible (after deforming) for any self shrinker, even of ``high'' entropy; topologically this seems to give no information not already known by the classical theory of surfaces. A finiteness theorem for self shrinking surfaces might be possible; one would have to find hypotheses in which one could gain uniform control on the noncollapsing constant. 
$\medskip$

In \cite{BHH} and \cite{BHH1}, Buzano, Haslhofer, and Hershkovits prove connectedness results for the moduli space of $2$-convex spheres in $\R^N$; by using the mean cuvature flow with surgery for surfaces of low entropy one sees that these results can be extended to show that the space of mean convex low entropy and spheres and tori are connected via mean convex paths. It would be interesting then to understand if honest connectedness results are possible in the low entropy mean convex case as well. 
$\medskip$

Relatedly, we remark it is easy to see that low entropy, mean convex hypersurfaces and 2-convex hypersurfaces are different subsets of the space of all embedded hypersurfaces. For example, by varying a small piece of a low entropy, mean convex hypersurface near a distance maximizing point from the origin, which must be locally convex, by ``flattening'' one can easily find a nearby hypersurface that is mean convex but not two convex, and the nearby hypersurface can be constructed to have low entropy as well since the entropy functionals concentrated near the maximal point where the variation is done will be approximately 1 throughout the isotopy by making the support of the variation small. So there are low entropy, mean convex hypersurfaces that aren't 2-convex; on the other hand, there are also 2-convex hypersurfaces that aren't low entropy. One way to construct such a 2-convex hypersurface that is not of low entropy is to consider a tubular neighborhood of a curve winding around the shrinking $S^n$: it is easy to construct a sequence of curves so that the area of these neighborhoods approaches twice the area of the shrinking sphere and so that the neighborhoods are 2 convex. The Gaussian area ($F_{0,1}$ functional) of these examples will be about twice that of $\Lambda_n$ then and hence bigger than 2 (since we always have $\Lambda_n > \sqrt{2}$) so not low entropy.
$\medskip$

\begin{appendix}

\section{the need for a carefully designed cap}

In this appendix we illustrate the need for a carefully designed cap by showing that the entropy of a capped off half cylinder  will be strictly greater than that of a cylinder; in particular the potential increase of the $F$ functionals near the surgery is a real issue and that care is warranted in understanding how much the entropy will increase. 
$\medskip$

In the following, we consider the toy surgery of capping a half cylinder off with a hemisphere. We only consider $F$ funcntionals for which the integral involved is particularly symmetrical; one quickly sees that calculating the $F$ functionals directly for a cap (even in the toy case, where the capped off cylinder can be explicitly parameterized) is onerous; thats why in the above we opt instead to let the mean curvature flow create (in a sense) the cap for us. 
$\medskip$

That being said, consider the diagram of the cap below, where $C$ is the center point of the hemisphere of radius 1, $a$ is the distance of the point $x_0$ from $C$ which lies in the center of the core of the cylinder, and the whole cap is denoted $\Sigma$:
\begin{center}
$\includegraphics[scale = .5]{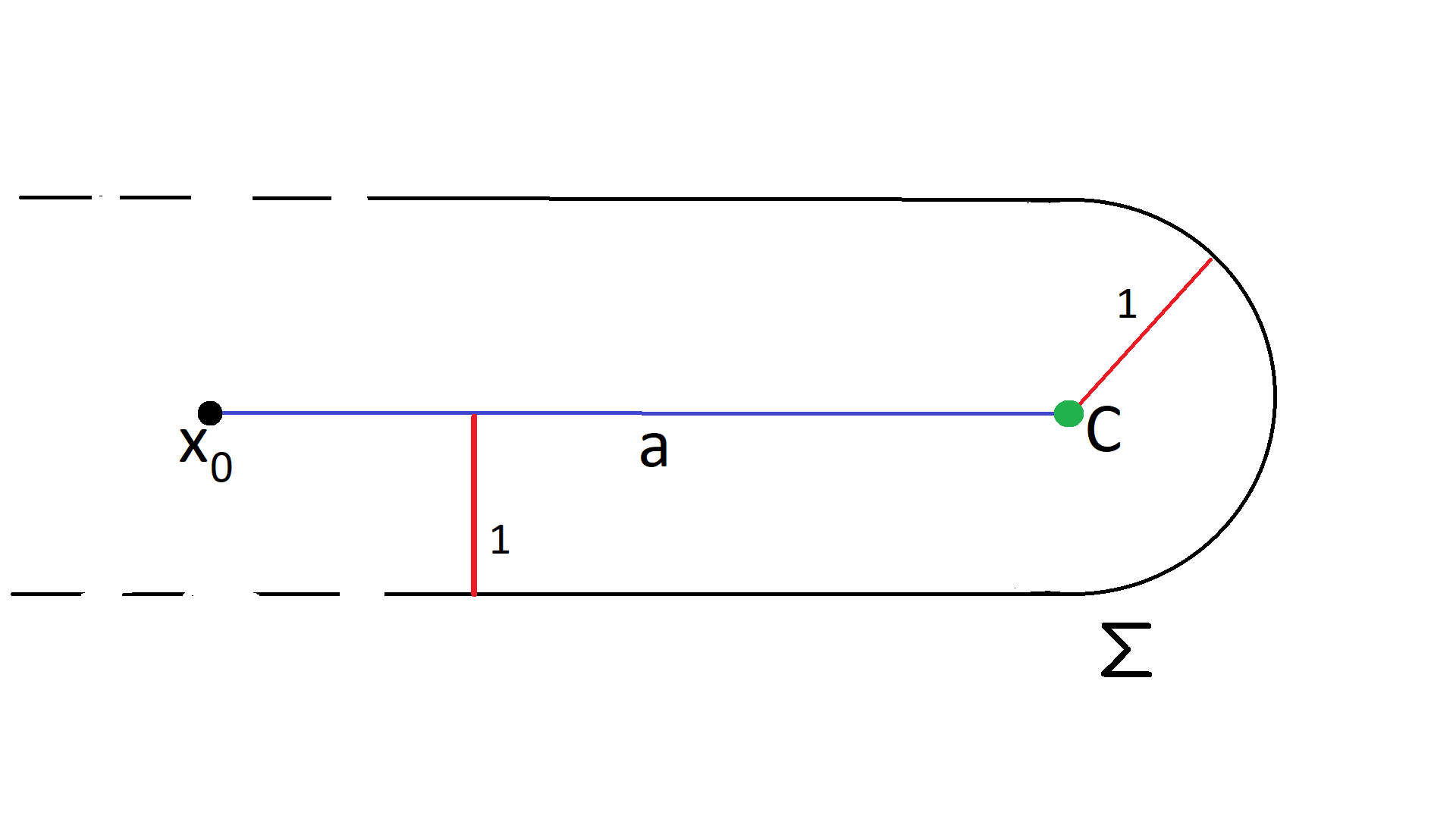}$
\end{center}
Now, on one hand we can see by taking the derivative of $F_{x_0(a),r}$ in $a$ that the $F$ functionals are a decreasing function of $a$; writing $F_{x_0, r}$ in terms of $a$:
\begin{equation}
\begin{split}
&F_{x_0, r}(\Sigma) = \int_{-a}^{-\infty} \frac{1}{4 \pi r} e^{\frac{-|x^2 + 1|}{4r}} 2 \pi dx +  \int_0^1 \frac{1}{4\pi r}e^{\frac{-|1 - x^2 + (a + x)^2|}{4r}} 2\pi \sqrt{1 - x^2} \frac{1}{\sqrt{1-x^2}} dx \\
& = \int_{-a}^{-\infty} \frac{1}{2r} e^{\frac{-|x^2 + 1|}{4r}} dx + \int_0^1 \frac{1}{2r} e^{\frac{-(a^2 + 2ax + 1)}{4r}} dx
\end{split}
\end{equation}
The derivative in $a$ then is given by:
\begin{equation}
\begin{split}
& \frac{d}{da} F_{x_0(a), r}(\Sigma) = -\frac{1}{2r} e^{\frac{-(a^2 + 1)}{4r}} + \frac{d}{da} \int_0^1e^{\frac{-(a^2 + 2ax + 1)}{4r}} dx \\
& =  \frac{-1}{2r} e^{\frac{-(a^2 + 1)}{4r}}  + \frac{d}{da}(\frac{e^{\frac{-a^2 -1}{4r}} - e^{\frac{-(a + 1)^2}{4r}}a}{a}) \\
& =  \frac{-1}{2r} e^{\frac{-(a^2 + 1)}{4r}} + \frac{e^{\frac{-(a + 1)^2}{4r}} \frac{a^2 + 2r + a}{2r} - e^{\frac{-a^2 -1}{4r}} \frac{a^2 + 2r}{2r}}{a^2} \\
& = \frac{-e^{\frac{-a^2 - 1}{4r}}}{a^2} + e^{\frac{-(a + 1)^2}{4r}} \frac{a^2 + 2r + a}{2a^2 r} \\
\end{split}
\end{equation}
We see for fixed choice of $r$ that for $a$ big enough, this is nonpositive. On the other hand, we see that as $a$ tends to infinity for a fixed scale $r$ that the $F$ functionals must tend to that of the cylinder; implying the $F$ functionals for $a$ finite must be strictly greater than that of the round cylinder and hence the entropy of the postsurgery domain must be too. 

\end{appendix}

\end{document}